\documentclass[prd,eqsecnum,showpacs,nofootinbib]{revtex4}
\usepackage{amsmath}
\usepackage{bm}
\usepackage{amssymb}
\usepackage{graphicx}
\usepackage{epstopdf}

\begin{document}
\title{Equivalence of Simplicial Ricci Flow and Hamilton's Ricci Flow\\ for 3D Neckpinch Geometries}
\author{Warner A. Miller$^{1,2}$, Paul M. Alsing$^{3}$, Matthew Corne$^{3}$ \& Shannon Ray$^{1}$}
\affiliation{
$^1$ Department of Physics, Florida Atlantic University, Boca Raton, FL 33431\\
$^2$ Department of Mathematics, Harvard University, Cambridge MA 02138 \\
$^3$ Air Force Research Laboratory, Information Directorate, Rome, NY 13441
}

\date{\today}
\begin{abstract} 
Hamilton's Ricci flow (RF) equations were recently expressed in terms of the edge lengths of a $d$-dimensional piecewise linear (PL) simplicial geometry, for $d\ge 2$.  The structure of the simplicial Ricci flow (SRF) equations are dimensionally agnostic.  These SRF equations were tested numerically and analytically in 3D for simple models and reproduced qualitatively the solution of continuum RF equations including a Type-1 neckpinch singularity.   Here we examine a continuum limit of the SRF equations for 3D neck pinch geometries with an arbitrary radial profile. We show that the SRF equations converge to the corresponding  continuum RF equations as reported by Angenent and Knopf.   
\end{abstract}
\maketitle

\section{Exploring Simplicial Ricci Flow in 3D}

Hamilton's Ricci flow (RF) continues to yield new insights into problems in pure and applied mathematics and proves to be a useful tool across a broad spectrum of engineering fields \cite{Cao:2003,Chow:2004,Chow:2006,Chow:2007}.  Here the time evolution of the metric is proportional to the Ricci tensor, 
\begin{equation}\label{eq:RF}
g^{ac}\, \dot g_{cb} = -2\, Rc^a{}_{b}.
\end{equation}
Hamilton showed that this yields a forced diffusion equation for the curvature, e.g. the scalar curvature evolves as 
\begin{equation}\label{eq:fdeqn}
\dot R = \triangle R + 2 R^2.
\end{equation}
The bulk of the applications of this curvature flow have utilized the numerical evolution of piecewise-flat simplicial 2-surfaces \cite{Gu:2012,Chow:2003}.  It is a widely accepted verity  in computational science that a geometry with complex topology is most naturally represented in a coordinate-free way by an unstructured mesh.  This is apparent in the engineering applications utilizing  finite-volume \cite{Peiro:2005} and finite-element \cite{Humphries:1997} algorithms, and is equally true in physics within the field of general relativity through Regge calculus \cite{Regge:1961,Gentle:1998}, and in electrodynamics by discrete exterior calculus  \cite{Hirani:2005}.  One expects a wealth of exciting new applications for discrete formulations of Ricci flow in 3 and higher dimensions.  Our presumption is  based on the uniformization theorem in 2-dimensions,  and since the topological taxonomy is richer in higher dimensions \cite{Thurston:1997}.   We envision that these higher-dimensional applications will involve geometries with complex topology and geometry.  Preliminary work in this direction is already underway  \cite{Glickenstein:2011,G:2005,Ge:2013,Forman:2003}.   In accord, we recently introduced a discrete RF approach for three and higher dimensions  that we refer to as simplicial Ricci flow (SRF) \cite{Miller:2013}.  SRF is founded on Regge calculus and upon the mathematical foundations of Alexandrov \cite{Regge:1961, Alexandrov:1950,AMM:2011,McDonald:2012}.  Here the SRF equations are similar to their continuum counterpart.  They are naturally defined on a $d$-dimensional simplicial geometry as a proportionality between the time rate of change of the circumcentric dual edges, $\lambda_i$, and the  simplicial Ricci tensor associated to these dual edges, 
\begin{equation}\label{eq:ref}
\dot \lambda_i /\lambda_i= -Rc_\lambda.
\end{equation}

It is the aim of this paper to explicitly show that these SRF equations in 3-dimensions for a geometry with axial symmetry converge to the continuum Hamilton RF equations for an arbitrarily refined mesh.
We examine a continuum limit of the SRF equations for 3D neck pinch geometries. We show that the SRF equations converge to the corresponding  continuum RF equations  as reported by  Angenent and Knopf \cite{Knopf:2004,Knopf:2011}.   In particular, we examine a piecewise flat axisymmetric geometry with 3-sphere topology, $S^3=S^2 \times \left[0,1\right]$.  Our lattice 3-geometry is tiled primarily with triangle-based frustum polyhedra with two simplicial polyhedra ``end caps.'' This axisymmetric lattice 3-geometry is characterized by edges of two kinds, the axial edges, $a_i$ and the cross-sectional sphere edges, $s_i$.  In recent work we examined a discrete model where each of the cross-sectional spheres was approximated by an icosahedron, of edge $s_i$ \cite{Alsing:2013}.  The accuracy of our model was limited by the relatively low spatial resolution of the icosahedron. The model considered here is built of infinitesimal isosceles-based frustums. The spatial resolution, $\ell$, of each of the spherical cross-sectional polyhedra with radius $\rho$ is driven by a single infinitesimal scale $\ell=\rho \xi \rightarrow 0$, and we assume an arbitrarily large number of cross-sectional spheres. The spherical edges $s_i$ and axial edges $a_i$ are infinitesimal.   The axisymmetry of this model allows us to consider only two SRF equations, one associated to each axial edge, $a_i$, the other with a spherical edge, $s_i$. In the limit there will be an ever increasing number of frustum blocks in our  model. We show here that we recover the exact RF equations in the continuum limit.  

The foundation of this work is the analysis of Angenent and Knopf on the neck pinching singularity analysis of  RF on a class of axisymmetric double-lobe-shaped geometries of Fig.~ \ref{fig:dumbbell}.  The metrics they considered were warped product metrics on $I \times S^2$, 
\begin{equation}\label{eq:wpm}
g = da^2 + \rho^2 d\Omega^2,
\end{equation}
where $I \in R$ is an open interval, $ d\Omega^2 = d\theta^2+\sin^2{(\theta)} d\phi^2,$ is the usual metric of the unit 2-sphere, $a$ is the proper axial distance away from the equator, and $\rho(a)$ is the cylindrical radial profile of the axi-symmetric geometry, i.e. $\rho(a)$ is the radius of the cross-sectional 2-sphere at an axial distance $a$ away from the equator. The vector-valued one-form Ricci tensor is,
\begin{equation}
Rc = e_a \left( -2 \frac{\rho''}{\rho} \right) \omega^a + e_\theta \left( \frac{1}{\rho^2}-\frac{\rho''}{\rho} - \left(\frac{\rho'}{\rho}\right)^2\right) \omega^\theta + 
e_\phi \left( \frac{1}{\rho^2}-\frac{\rho''}{\rho} - \left(\frac{\rho'}{\rho}\right)^2\right) \omega^\phi,
\end{equation}
where the primes refer to partial derivatives with respect to the axial distance, $a$, and the dots refer to derivatives with respect to time, $t$.
In particular, we show in this manuscript that the SRF equations, under a suitable mesh refinement,  converge to the continuum Hamilton's RF equations of Knopf and Angenent, 
\begin{eqnarray}
\label{eqn:crfa}
\frac{\dot{a}}{a} & =  -Rc^a{}_a  & = 2\, \frac{\rho''}{\rho},\\
\label{eqn:crfrho}
\frac{\dot{\rho}}{\rho} & =  -Rc^\theta{}_\theta & = \rho''/\rho + \left(\rho'/\rho\right)^2 - 1/\rho^2 .
\end{eqnarray}
Armed with this result, all continuum theorems and corollaries apply equally to the SRF equations. They are equivalent in the continuum limit we use. Therefore, we can confidently say that the SRF equations are consistent with the following theorems of Angenent and Knopf for non-degenerate neck pinches:
\begin{enumerate}
\item If the scalar curvature is everywhere positive, $R\ge0$,  then the radius of the waist ($a_{min}=\rho(0)$)  is bounded,
\begin{equation}
(T-t) \le a_{min}^2 \le 2(T-t),
\end{equation}
where $T$ is the finite time that a neck pinch occurs.
\item As a consequence the neck pinch singularity occurs at or before, $T = a_{min}^2$.
\item The height of the two lobes is bounded from below, and under suitable conditions the neck will pinch off before the lobes will collapse. 
\item The neck approaches a cylindrical-type singularity.
\end{enumerate}
The result of this paper shows that any continuum RF theorem or curvature bound  for this class of geometry will apply equally well to the SRF evolution with this geometry -- the SRF equations are  equivalent to the RF equations in a continuum limit.  Our result applies to both degenerate and non-degenerate neck pinch singularity formation. The proof of the  convergence of SRF to continuum RF is done here  explicitly and algebraically.    While this does not prove the equivalence between Hamilton's RF and SRF equations for any geometry and for any dimension; nevertheless, we conjecture this is true.  The work here supports the definition (Def.~1) of the SRF equations introduced recently in \cite{Miller:2013}.

\section{A Lattice Approximation of the Angenent-Knopf Neckpinch Geometry}
\label{sec:2}

\begin{figure}
\includegraphics[width=7.5cm]{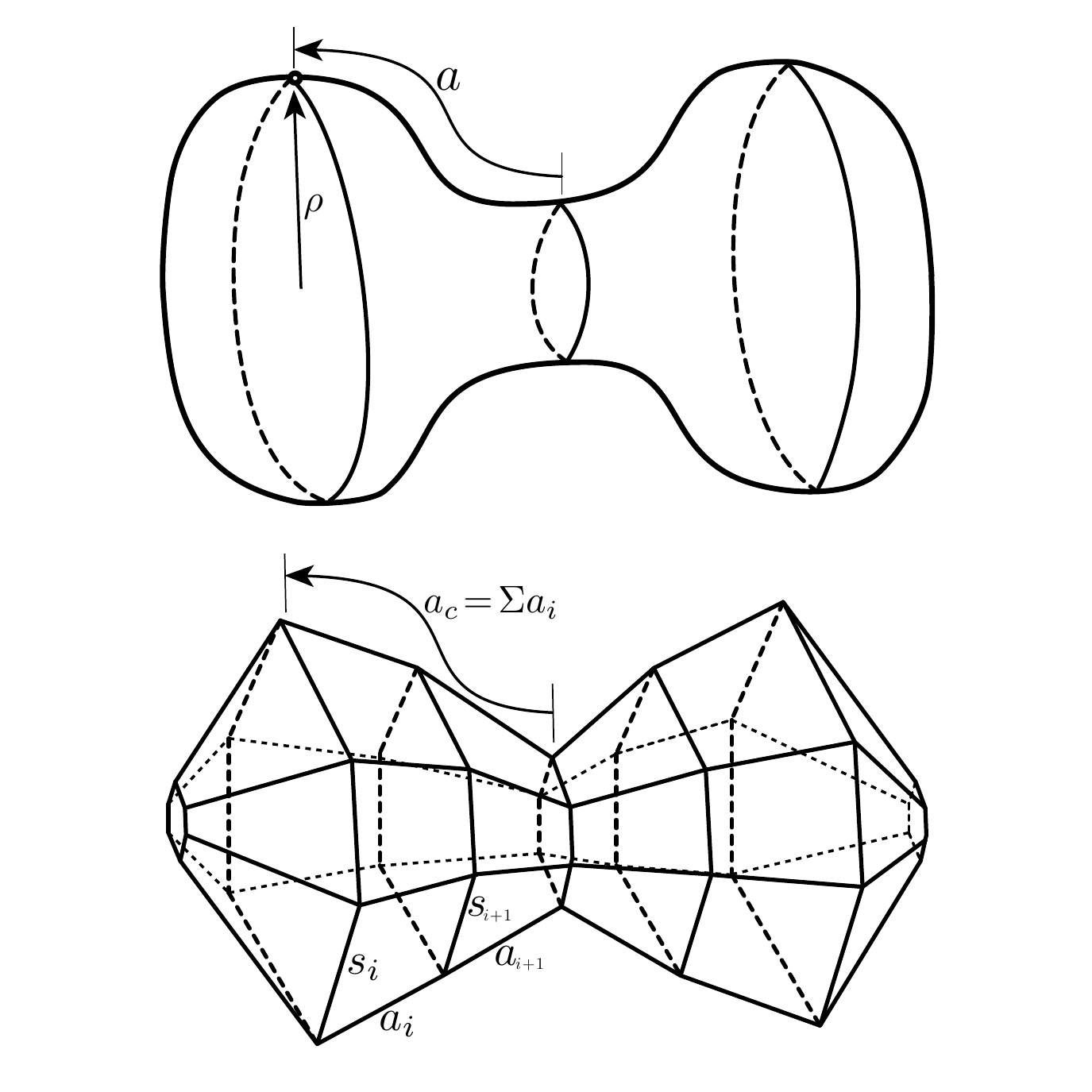}
\caption{A 2-dimensional rendering of the simplicial double-lobed model.  For visualization purposes, we have suppressed one of the azimuthal angles for each of the 2-sphere cross sections. In two dimensions the piecewise linear (PL) surface tiles are trapezoids, where in 3-dimensions they are frustum polyhedra (see Fig.~\ref{fig:frustum} in the Appendix).  In this paper we take the limit of an ever more finely discretized sphere, and with an ever increasing number of spherical cross sections. In this limit we show that the continuum Hamiltonian RF equations are recovered from the SRF equations.}
\label{fig:dumbbell}
\end{figure}

For the purpose examining Type-1, or Type-2  neck pinch behavior of the SRF equations, we have introduced a PL lattice geometry sharing the qualitative features of the Angenent and Knopf initial data \cite{Knopf:2004} as illustrated in Fig.~\ref{fig:dumbbell}.  The continuum cross sections of this geometry in planes perpendicular to the symmetry axis are 2-sphere surfaces. We impose no mirror symmetry in the radial profile, $\rho=\rho(a)$.  The surface and metric can be parameterized by two coordinates, $a$ and $\rho(a)$. Here a given point on the surface is identified by its proper ``axial" distance of $a$ from the equator (or any other point we so choose) , and the radius, $\rho=\rho(a)$, of the cross-sectional sphere on which the point lies. The continuum warped-product metric of this surface as introduced by Angenent and Knopf \cite{Knopf:2004},
\begin{equation}\label{eq:metric}
 g_{ij} = da^2+\rho(a)^2 d\Omega^2,
 \end{equation}
 where
 \begin{equation}
 d\Omega^2 = d\theta^2+\sin^2{(\theta)} d\phi^2,
 \end{equation}
 is the usual spherical line element.
The initial data is determined by a radial profile function at $t=0$ for the double-lobed geometry, and amounts to specifying a function relating the cylindrical radius, $\rho$,  to a scaled proper axial distance along the double-lobed geometry away from an equator or neck, $a\in\{a_{min},a_{max}\}$,
\begin{equation}
s = \psi(a).
\end{equation} 
By way of an example, if the double-lobed geometry has no neck, and were just a sphere of radius, $R_0$, this initial  radial profile function is simply the cylindrical coordinate radius, 
\begin{equation}
\rho(a,t=0)=R_0\, cos(a).
\end{equation}
However, Angenent and Knopf's mirror-symmetric double-lobed geometry introduced a parabolic waist for their purpose so as to aid in their mathematical analysis of the neck singularity, 
\begin{equation}
\rho=\rho_{AK} = \left\{
\begin{array}{l l}
R_0 \cos{(a)}  & |a|\ge \frac{\pi}{4} \\
R_0 \sqrt{A + B a^2} & |a| < \frac{\pi}{4}
\end{array}
\right. ,
\end{equation}
where constant $A$  controls the degree of neck pinching in the initial-value data, and the constant $B$ is chosen so as to ensure continuity in the radial profile function at $a=\pm \pi/4$. In this paper, we explore arbitrary radial profiles.  

\begin{figure}
 \includegraphics[width=10.5cm]{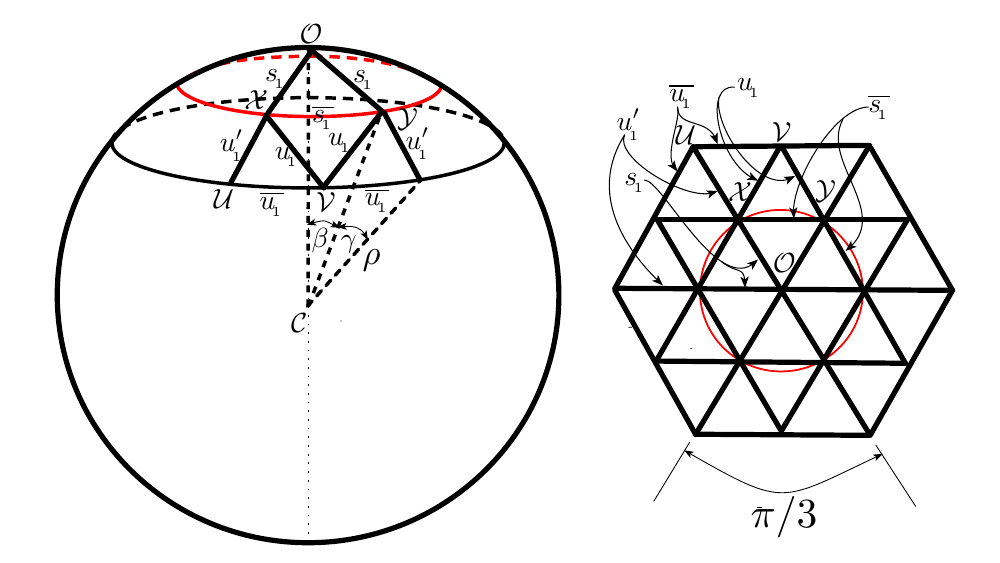}
\caption{The 2-dimensional triangulated region around the vertex,  ${\cal O}$ , on a spherical cross-section of the axisymmetric geometry.  In this manuscript we assume $\xi=\ell/r$ is a global infinitesimal for each of the, $N_a$, triangulated sphere cross sections in this model. }
\label{fig:deltatheta}
\end{figure}

\begin{figure}
 \includegraphics[width=10.5cm]{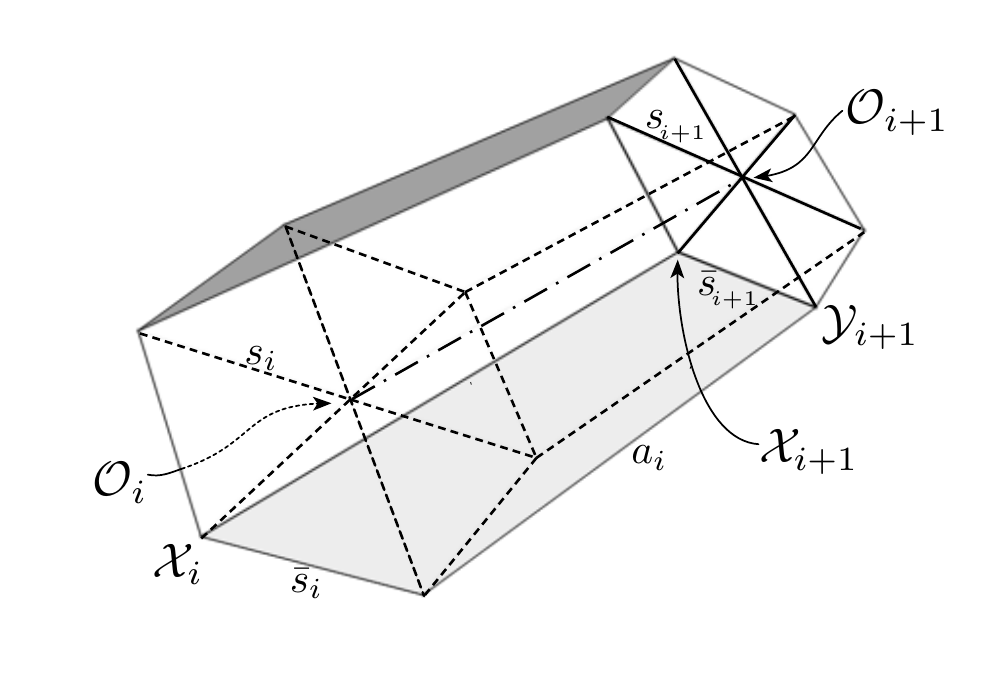}
\caption{The six identical isosceles frustum blocks sharing edge $a_i=\overline{{\cal O}_i {\cal O}_{i+1}}$.}
\label{fig:aiedge}
\end{figure}

We provide a simplicial approximation of  an axisymmetric warped-product geometry at time $t$ characterized by an arbitrary $C^2$ radial profile
\begin{equation}
\rho(a,t)=\rho(a)\ \forall \, a\in\{a_{min},a_{max}\}.
\end{equation}
We first identify an arbitrarily large number ($N_s \rightarrow \infty$) of nearly equal-spaced spherical cross sections. Next we examine one of these spheres, namely the $i'th$ cross-sectional sphere. There are many ways to approximate this by a polyhedron with an arbitrarily large amount of vertices, $N_{s} \rightarrow \infty$.  We utilize the symmetry of our model to concentrate only in the vicinity of a single point $\cal O$ on the sphere.  At this point we project a infinitesimal flat-space hexagonal lattice onto its surface, as shown in Fig.~\ref{fig:deltatheta}. Here we take the length of the isosceles triangles that we are projecting onto the sphere of radius $\rho$  to be arbitrarily small, $\ell_k \ll \rho$.  This yields an infinitesimal parameter in our model that we will drive to zero, 
\begin{equation} 
\xi := \frac{\ell}{\rho} \rightarrow 0.
\end{equation}
In order to construct the SRF equations at $\cal O$ it is necessary for us to extend the lattice radially one more level out from the first six equilateral triangles so that we can examine 18 additional equilateral triangles, each of edge length $\ell_k \ll \rho$, that we project onto the surface of the sphere as shown in the right-hand side of  Fig.~\ref{fig:deltatheta}.  The projected triangles will no longer be equilateral. In particular there will be two sets of 6 isosceles triangles, $\{\cal O, X, Y\}$ and $\{\cal V, X, Y\}$, as well as twelve triangles with three different edge lengths $\{\cal X,U,V\}$.  These 24 triangles are composed of combinations of six distinct edges,
\begin{eqnarray}
\label{eqn:s}
s  = &  \arccos{\left( \frac{\overrightarrow{\cal CO}}{|{\cal CO}|} \cdot \frac{\overrightarrow{\cal CX}}{|{\cal CX}|}\right) } &\approx  \rho\, \xi  \left( 1 - \frac{1}{3} \xi^2 + \frac{1}{5} \xi^4 + O\left[ \xi \right]^6\right),\\
\bar s  = &  \arccos{\left( \widehat{\cal CX}\cdot \widehat{\cal CY} \right) }  &\approx \rho\, \xi \left( 1 - \frac{11}{24} \xi^2 + \frac{203}{640} \xi^4 + O\left[ \xi \right]^6\right),\\
u  = &  \arccos{\left(\widehat{\cal CX}\cdot \widehat{\cal CV}\right) } &\approx \rho\, \xi \left( 1- \frac{35}{24} \xi^2 + \frac{1183}{640} \xi^54+ O\left[ \xi \right]^6\right),\\
\bar u  = &  \arccos{\left( \widehat{\cal CU}\cdot \widehat{\cal CV}\right) } &\approx  \rho\, \xi \left( 1 - \frac{11}{6} \xi^2 + \frac{203}{40} \xi^4 + O\left[ \xi \right]^6\right),\\
\label{eqn:up}
u'  = &  \arccos{\left(\widehat{\cal CX}\cdot \widehat{\cal CU}\right) } &\approx \rho\, \xi \left( 1 - \frac{7}{3} \xi^2 + \frac{31}{5} \xi^4 + O\left[ \xi \right]^6\right).
\end{eqnarray} 
We assume that all of the triangulated spherical polyhedral cross-sections in our model have the same lattice topology.  Furthermore, we assume that they are all congruent to each other under a suitable  global  scale factor.  Consider the $i$'th and $(i+1)$'st polyhedral sphere of radius $\rho_i$ and $\rho_{i+1}$; respectively. Each triangulated polyhedron has $N_0 \gg 1$ vertices, and therefore $N_1=3 N_0 - 6$ edges and $N_2=2 N_0 -4$ triangles.  We connect these two polyhedra together by connecting the $N_0$ pairs of corresponding pairs of vertices by $N_0$ identical axial edges, each of length $a_i$.  Each pair of corresponding triangles when connected by three $a_i$ edges will form a triangular-based frustum block (see Appendix, Fig.~\ref{fig:frustum}).  We require that the geometry interior to each of our frustum blocks is flat Euclidean 3-space.  Consequently, the 3-dimensional geometry between the two bounding spherical polyhedrons is tiled with $N_2$ frustum blocks, one for each of the triangles. In particular,  axial edge $\overline{{\cal O}_i{\cal O}_{i+1}} = a_i$ is the meeting place of  six identical isosceles triangle frustum blocks as illustrated in Fig.~\ref{fig:aiedge}.  This is not ordinarily the case for every axial edge, $a_i$, e.g.  axial edge $\overline{{\cal X}_i{\cal X}_{i+1}} = a_i$ is the meeting place of three distinct pairs of  frustum blocks, two pairs are isosceles frustum blocks the last pair is a general triangular-based frustum block.  
The axial symmetry of this model permits us to tile the geometry with the non-simplicial furstum blocks. Given the symmetry of our model the geometry of each frustum block is completely determined by its 9 edge lengths, in other words the symmetry endows each frustum block with rigidity.  This construction gives us an axisymmetric 3-cylinder geometry composed of triangular-based frustums.  We ``cap-off" this 3-cylinder by treating each of the two bounding polyhedrons of radius $\rho_1$ and $\rho_{N_a}$ as flat Euclidean tiles of our PL geometry.  Our lattice geometry becomes homeomorphic to a 3-sphere when we include the two polyhedral end caps.  It is composed of $N_2( N_a-1)$ frustum blocks and two triangulated polyhedral ``end caps.''   This can be seen in Fig.~\ref{fig:dumbbell}, albeit with one dimension suppressed. There the end caps are flat hexagons as opposed to triangulated polygons.   

A continuum limit of our lattice is achieved here by taking (1) $N_a \rightarrow \infty$, or equivalently $a_i \rightarrow 0$, and (2)   $N_s \rightarrow \infty$, or equivalently $\xi \rightarrow 0$.

 \section{Theorem: The SRF Equations for the Frustum Geometry are the Hamilton RF equations}
 \label{sec:theorem}
 
The recent definition of the SRF equations and the dual-edge SRF equations in \cite{Miller:2013} made use of elements from both the simplicial lattice geometry (in this case the frustum lattice geometry with end caps, $\cal F$), and  from its circumcentric dual lattice $\cal F^*$.  The SRF and dual-edge SRF equations for the $S^3$ frustum geometry, $\cal F$ are functions of the $(N_a - 1)$ axial edges, $a_i\ \forall i\in\{1,2\ldots,N_a-1\}$, as well as the $N_a$ radii $\rho_i\ \forall i\in\{1,2\ldots,N_a\}$.   The dual circumcentric lattice is composed of the dual axial edges, $\alpha_i\ \forall i\in\{1,2\ldots,N_a\}$, and the dual edges $\sigma_i\ \forall i\in\{1,2\ldots,N_a+1\}$. Each dual axial edge reaches from the circumcenter of one frustum block to the circumcenter of the adjacent frustum block sharing a common triangle face. The dual edge $\alpha$ pierces the triangle at the triangles circumcenter, and the edge is perpendicular to the triangle.  On the other hand, the dual spherical edge $\sigma$ reaches from the circumcenter of one frustum block to the circumcenter of an adjacent frustum block sharing a common trapezoidal face.  This dual edge $\sigma$ is perpendicular to the trapezoid and pierces the trapezoid at its circumcenter.  Therefore, there will be two dual axial edges, $\alpha$ and three dual spherical edges $\sigma$ emanating form the circumcenter of each of the frustum blocks.  While the circumcenter of each of the two polyhedral ``end caps'' will be the common meeting place of the $N_2$ dual axial edges.  

We proved recently that the dual-edge SRF equations are equivalent to the simplicial-edge SRF equations for this warped-product geometry in Appendix B of \cite{Alsing:2013}, 
\begin{equation}\label{eq:drrf-rrf}
\underbrace{
\left\{
\begin{array}{l}
\frac{\dot \sigma_i}{\sigma_i} = -Rc_{\sigma_i}\\
\frac{\dot \alpha_i}{\alpha_i} = -Rc_{\alpha_i}
\end{array}
\right\}
}_{dual-edge\ SRF\ equations}
\equiv \ \ \
\underbrace{
\left\{
\begin{array}{l}
\sum_{\lambda_j \in s^*_i}  \frac{\dot \lambda_j}{\lambda_j} \left(\frac{V_{s_i \lambda_j}}{V_{s_i}}\right) = -Rc_{s_i}\\
\sum_{\lambda_j \in a^*_i}  \frac{\dot \lambda_j}{\lambda_j} \left(\frac{V_{a_i \lambda_j}}{V_{a_i}}\right) = - Rc_{a_i}
\end{array}
\right\}
}_{simplicial-edge SRF\ equations},
\end{equation} 
where the definition of the volume-weighting factors are defined below in \ref{sec:vwf} .
Therefore, and for the purpose of this paper, it will suffice to prove the following theorem:\\
\\
\noindent
{\bf Theorem 1}\\
\noindent 
The two dual-edge SRF (Eq.~\ref{eq:ref}) equations at a vertex ${\cal O}_i$  in the frustum-based warped-product lattice geometry, $\cal F$,  converge to their continuum RF counterparts, Eqs.~\ref{eqn:crfa}-\ref{eqn:crfrho}, 
 \begin{eqnarray}
 \dot \sigma_i /\sigma_i = -Rc_{\sigma_i} & \Longrightarrow & \frac{\dot{\rho}}{\rho}  =  \rho''/\rho + \left(\rho'/\rho\right)^2 - 1/\rho^2,\\
 \dot \alpha_i /\alpha_i= -Rc_{\alpha_i} &  \Longrightarrow & \frac{\dot{a}}{a}  =   2\, \frac{\rho''}{\rho},
 \end{eqnarray} 
 in the limit when (1) $\xi \rightarrow 0$, and (2) $N_a \rightarrow \infty$. 
 
 \section{The Dual $\alpha_i$-Edge SRF Equation at the Point ${\cal O}_i$ and its Continuum Limit }
 \label{sec:alphaSRF}

\begin{figure}
\includegraphics[width=8.5cm]{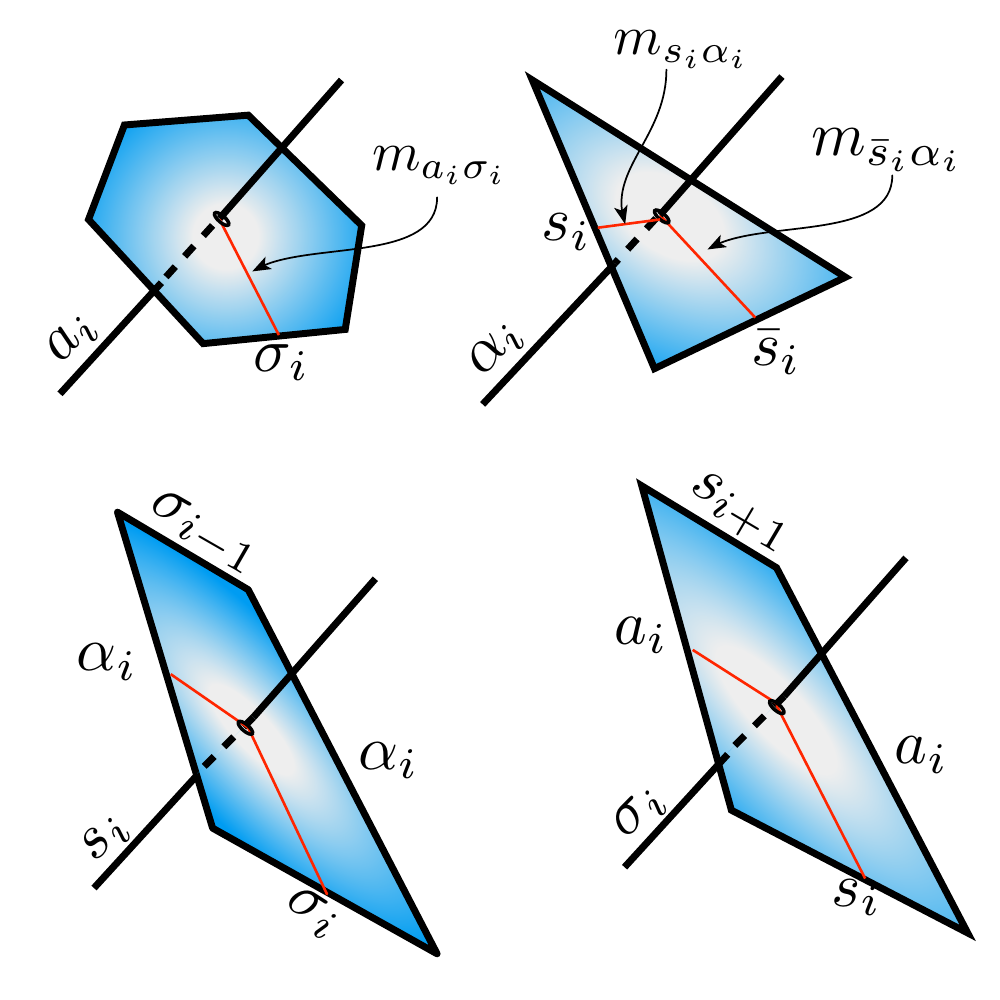}
\caption{In 3-dimensions each edge of the triangular-based frustum lattice is dual to an area of the dual circumcentric lattice. This is shown in the first column for edges $a_i$  and $s_i$. Conversely, each dual edge in the circumcentric dual  lattice is dual to an area of the triangular-based frustum lattice, as illustrated in the second column for edges $\alpha_i$  and $\sigma_i$. }
\label{fig:dualareas}
\end{figure}

In this section we examine the dual-edge SRF equation associated to the axial edge, $\alpha_i\in {\cal F}^*$  as displayed in the lower left-hand side of Eq.~\ref{eq:drrf-rrf}.  Edge $\alpha_i$  is dual to triangle $\bigtriangleup_{{\cal O}_i {\cal X}_i {\cal Y}_i} \in {\cal F}$.  The edge $\alpha$ and triangle $\bigtriangleup_{{\cal O}_i {\cal X}_i {\cal Y}_i}$ are illustrated in the upper right hand side of              Fig.~\ref{fig:dualareas}.  We also use Eqs.~\ref{eqn:s}-\ref{eqn:up} to examine the continuum limit of this equation. We use the definition of the dual-edge SRF equations as introduced recently in \cite{Miller:2013}, and consequently we find
\begin{equation}
\label{eqn:alphaequation}
\frac{\dot \alpha_i}{\alpha_i} = -Rc_{\alpha_i }  = -\sum_{\ell_j\in\alpha^*_i} 2 \frac{\epsilon_{\ell_j}}{\ell^*_j} \left(  \frac{V_{\alpha_i \ell_j}}{V_{\alpha_i}}\right) = -4 \frac{\epsilon_{s_i}}{s^*_i} \left(  \frac{V_{\alpha_i s_i}}{V_{\alpha_i}}\right) - 2 \frac{\epsilon_{\bar s_i}}{\bar s^*_i} \left(  \frac{V_{\alpha_i \bar s_i}}{V_{\alpha_i}}\right). 
\end{equation}
We use the results of the Appendix and the definitions in \cite{Miller:2013} to examine each term of this equation and to examine the Taylor series expansion to obtain the continuum limit expression. We do this in three steps by examining (1) the two circumcentric volume weighting factors on the right-hand side of the dual SRF equation, (2) the two Gaussian curvature expressions also on the right-hand side of Eq.~\ref{eqn:alphaequation}, and finally (3) the time derivative on the left-hand side of this equation. 

\subsection{The Dual Edge $\alpha_i$ to Triangle $\bigtriangleup_i := \overline{{\cal O}_i {\cal X}_i {\cal Y}_i }$ }

The dual edge $\alpha_i$ reaches from the circumcenter  ${\cal C}_{3_{i-1}}$ of one frustum block  ${\cal F}_{i-1}$ to the circumcenter  ${\cal C}_{3_{i}}$  of an adjacent frustum block ${\cal F}_{i}$.   This dual edge is perpendicular to triangle  $\bigtriangleup_i := \overline{{\cal O}_i {\cal X}_i {\cal Y}_i }$ common to ${\cal F}_{i-1}$ and ${\cal F}_{i}$. From \ref{subsec:dualsegments}, Eqs.~\ref{eqn:halfalphaibase}-\ref{eqn:halfalphaitop}, we find
\begin{align}
\label{eqn:alphai}
\alpha_i & = \alpha_{top_{i-1}} + \alpha_{base_i} \\
& = 
\frac{
8 a_i^2 \bigtriangleup_i^2 + \bar s_i^2 {s'}_i^2 s_i \left(s_{i-1}-s_{i}\right)
}{
4 \bigtriangleup_i \sqrt{16 a_i^2 \bigtriangleup_i^2 - \bar s_i^2 {s'}_i^2 \left(s_i-s_{i-1}\right)^2}
}
+
\frac{
8 a_i^2 \bigtriangleup_i^2 + \bar s_i^2 {s'}_i^2 s_i \left(s_{i+1}-s_{i}\right)
}{
4 \bigtriangleup_i \sqrt{16 a_i^2 \bigtriangleup_i^2 - \bar s_i^2 {s'}_i^2 \left(s_i-s_{i+1}\right)^2}
} \\
& \approx  \frac{a_{i-1}+a_i}{2} \left( 1 + \frac{a_i \left( \rho_{i-1}-\rho_i \right) \left(\rho_{i-1}+3 \rho_i\right) + a_{i-1} \left( \rho_{i}-\rho_{i+1} \right) \left(\rho_{i+1}+3 \rho_i\right)}{6 a_{i-1} a_i \left( a_{i-1}+a_i\right) }\ \xi^2 + O[\xi]^4 \right). 
\end{align}
We will display the next higher order terms in our expansions.   They may be useful for numerical applications or to study singularity formation analytically. 

\subsection{Circumcentric Volume Weighting Factors}
\label{sec:vwf}

The hybrid volume $V_{\alpha_i}$ is the volume of the hybrid polyhedron formed by the product of dual edge $\alpha_i\in {\cal F}^*$ and triangle $\bigtriangleup_{{\cal O}_i {\cal X}_i {\cal Y}_i} \in {\cal F}$ and  is the sum of there reduced hybrid tetrahedra,
\begin{equation}
V_{\alpha_i} = 2 V_{\alpha_i s_i} + V_{\alpha_i \bar s_i} = \frac{1}{3} \alpha_i\, s_i \,m_{\alpha_i s_i}  + \frac{1}{6} \alpha_i\, \bar s_i\, m_{\alpha_i \bar s_i}.
\end{equation}
The moment arms, $m_{\alpha_i s_i}$ and $m_{\alpha_i \bar s_i}$  are displayed in the upper right-hand part of Fig.~\ref{fig:dualareas}, and are calculated in Eqs.~\ref{eqn:msialphai}-\ref{eqn:msibaralphai}. The fractional volumes and their series expansion in terms of $\xi$ are, 
 \begin{eqnarray}
 \label{eqn:Valphaisi}
 \frac{V_{\alpha_i s_i}}{V_{\alpha_i}} & = \frac{ s_i m_{\alpha_i s_i} }{ 2\, \alpha^*_i }  & \approx  \frac{1}{3} \left(1-\frac{1}{12} \xi^2+ \frac{1}{16} \xi^4 + O\left[\xi\right]^6\right)\\
 \label{eqn:Valphaisbari}
 \frac{V_{\alpha_i \bar s_i}}{V_{\alpha_i}} &  = \frac{ \bar s_i m_{\alpha_i \bar s_i} }{ 2\, \alpha^*_i }  & \approx  \frac{1}{3} \left(1-\frac{1}{6} \xi^2+ \frac{1}{8} \xi^4 + O\left[\xi\right]^6\right).
 \end{eqnarray}

\subsection{Gaussian Curvature}

\begin{figure}
 \includegraphics[width=10.5cm]{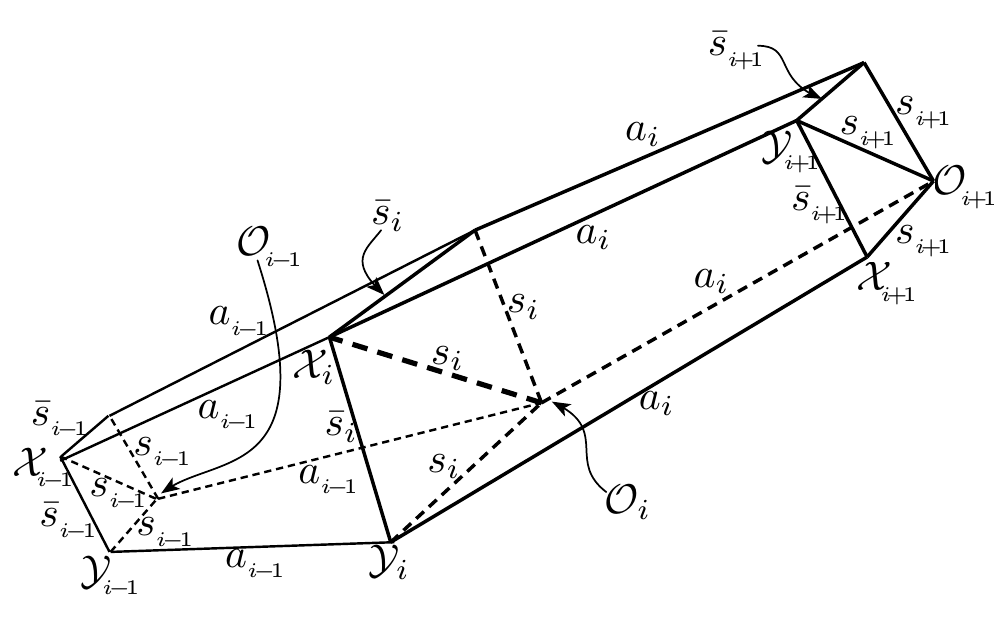}
\caption{The two pairs of isosceles frustum blocks sharing edge $s_i=\overline{{\cal O}_i {\cal X}_{i}}$.}
\label{fig:siedge}
\end{figure}

In \cite{Miller:2013} we defined the Gaussian curvature of edge $s_i$ to be the deficit angle, $\epsilon_{s_i}$, and distributed uniformly over the circumcentric dual area $s^*_i$, 
\begin{equation}
K_{s_i} := \frac{\epsilon_{s_i}}{s^*_i}.
\end{equation}
The edge, $s_i$, is common to four isosceles frustum blocks as shown in Fig.~\ref{fig:siedge}.
There are two pairs of dihedral angles of the frustum blocks at edge $s_i$, and we can use Eq.~\ref{eqn:thetasi} and Eq.~\ref{eqn:thetasip1} to define the deficit angle, 
\begin{align}
\epsilon_{s_i} & = 2\pi -2 \arccos{\left( \frac{\bar s_i \left(s_i-s_{i+1}\right)}{\sqrt{4 s_i^2-\bar s_i^2} \sqrt{4 a_i^2-\left(s_i-s_{i+1} \right)^2}}\right)}
-2 \arccos{\left( \frac{\bar s_{i} \left(s_i+s_{i-1}\right)-2 \bar s_{i-1} s_i}{\sqrt{4 s_i^2-\bar s_i^2} \sqrt{4 a_{i-1}^2-\left(s_i-s_{i-1} \right)^2}}\right)} \\
& \approx  \left(   \frac{   a_i\left(\rho_i-\rho_{i-1}\right)+a_{i-1} \left(\rho_i-\rho_{i+1} \right)    }{    \sqrt{3}\, a_i\, a_{i-1}    }\, \xi \right) \times \\ 
& \ \ \  \left( 1 + \frac{5 a_i^2\left(\rho_{i-1}-\rho_i\right)^2+a_{i-1}^2 \left(5 \left(\rho_i-\rho_{i+1}\right)^2-18 a_i^2\right) + 5 a_{i-1} a_i \left(\rho_{i-1}-\rho_i\right)\left(\rho_i-\rho_{i+1}\right)}{36\, a_i^2\,  a_{i-1}^2 } \xi^2 + O\left[\xi\right]^4 \right)  .
\end{align}
The dual area $s^*_i$ is a trapezoid and is shown in the lower left diagram of Fig.~\ref{fig:dualareas} and is the sum of four isosceles triangles, 
\begin{align}
s_i^* & = \frac{1}{2} \sigma_i m_{s_i \sigma_i} +  \frac{1}{2} \sigma_{i-1} m_{s_i \sigma_{i-1}} + \alpha_i m_{s_i \alpha_i} \\
& \approx \left(   \frac{   a_{i-1}\left(\rho_{i-1}+3 \rho_{i}\right)+a_{i} \left(3 \rho_i+\rho_{i+1} \right)    }{    8\sqrt{3}   }\, \xi \right) \times \\
& \ \ \ \left( 
1+ 
\frac{
\begin{array}{c}
a_{i} \left( \rho_{i-1}-\rho_{i}\right)\left( 5 \rho_{i-1}^2+10 \rho_{i-1} \rho_{i} +13 \rho_{i}^2\right) -6 a_{i-1}^2 a_{i}\left(\rho_{i-1}+3 \rho_{i}\right) +\\ 
a_{i-1} \left( -13 \rho_{i}^2+3 \rho_{i}^2 \rho_{i+1}+5 \rho_i \rho_{i+1}^2 + 5 \rho_{i+1}^3 -6 a_i^2\left( 3 \rho_i+\rho_{i+1}\right)\right) 
\end{array}
}{
12 a_i a_{i+1} \left( a_{i-1}\left(\rho_{i-1}+3 \rho_i\right) + a_i \left( 3 \rho_i + \rho_{i+1}\right)\right)
}
\xi^2
+ O\left[\xi\right]^4
\right).
\end{align}
The Gaussian curvature of edge $s_i$ in our lattice is expressed as, 
\begin{align}
\label{eqn:Ks}
K_{s_i} & := \frac{\epsilon_{s_i}}{s^*_i} \approx 
\frac{
8\left( a_i\left(\rho_i-\rho_{i-1}\right) - a_{i-1} \left( \rho_{i+1}-\rho_i \right) \right)
}{
a_i a_{i-1} \left( a_i \left(3 \rho_i +\rho_{i+1}\right) + a_{i-1} \left( 3 \rho_i + \rho_{i-1}\right) \right)
} \times 
\\
& \left(1+ 
\frac{
\begin{array}{c}
   (a_{i-1} (\rho_i-\rho_{i+1})+a_i (\rho_i-\rho_{i-1})) \left(5 a_{i-1}^2 (\rho_{i-1}+3
   \rho_i) (\rho_i-\rho_{i+1})+\right.\\
   \left. 2 a_{i-1} a_i \left(5 \rho_{i-1}^2+10 \rho_{i-1} \rho_i+12
   \rho_i^2+10 \rho_i \rho_{i+1}+5 \rho_{i+1}^2\right)-5 a_i^2 (\rho_{i-1}-\rho_i) (3
   \rho_i+\rho_{i+1})\right)
   \end{array}
   }{
   36 a_{i-1}^2 a_i^2 (a_{i-1} (\rho_{i-1}+3 \rho_i)+a_i (3\rho_i+\rho_{i+1}))
   }
   \, \xi ^2
   + O\left[\xi\right]^4 \right) .
\end{align}
Similarly we can calculate the Gaussian curvature curvature associated with the hinge edge $\bar s_i$ where
\begin{equation}
K_{\bar s_i} := \frac{\epsilon_{\bar s_i}}{\bar s*_i}.
\end{equation}
The dual area and dihedral angles now extend into the next band of triangles away from the pole ${\cal O}_i$ as shown in Fig.~\ref{fig:deltatheta}.
Using the expressions in the Appendix and the series expansion for the edges (Eq.~\ref{eqn:s}) the deficit angle for $\bar s_i$ is 
\begin{align}
\epsilon_{\bar s_i} & = 2\pi -
\arccos{\left(  \frac{  \left(s_{i}-s_{i+1}\right) \left(2 s_i^2 -\bar s_i^2 \right)  }{  s_i  \sqrt{4 s_i^2-\bar s_i^2} \sqrt{4 a_i^2-\left(\bar s_i-\bar s_{i+1} \right)^2}}\right)} - 
\arccos{\left(  \frac{  \left(u_{i}-u_{i+1}\right) \left(2 u_i^2 -\bar s_i^2 \right)  }{  u_i  \sqrt{4 u_i^2-\bar s_i^2} \sqrt{4 a_i^2-\left(\bar s_i-\bar s_{i+1} \right)^2}}\right)} - \\
&\arccos{\left(  \frac{  -\left(s_{i}-s_{i-1}\right) \left(2 s_i^2 -\bar s_i^2 \right)  }{  s_i  \sqrt{4 s_i^2-\bar s_i^2} \sqrt{4 a_i^2-\left(\bar s_i-\bar s_{i-1} \right)^2}}\right)} -
\arccos{\left(  \frac{  -\left(u_{i}-u_{i-1}\right) \left(2 u_i^2 -\bar s_i^2 \right)  }{  u_i  \sqrt{4 u_i^2-\bar s_i^2} \sqrt{4 a_i^2-\left(\bar s_i-\bar s_{i-1} \right)^2}}\right)}  \\
& \approx \frac{\left(\rho_{i-1}-\rho_{i+1}\right) \xi}{\sqrt{3}\, a_{i}} \left(1 +\frac{\left( 10\left(\rho_{i-1}^2+3 \rho_i^2 - 3 \rho_i \rho_{i+1} + \rho_{i+1}^2-3 \rho_{i-1} \rho_i +\rho_{i-1}\rho_{i+1}\right) -117 a_i^2  \right) \xi^2}{72 a_i} +O\left[\xi\right]^4\right), 
\end{align}
and the dual area $\bar s^*_i$ is the sum of four kite areas whose series in $\xi$ is
\begin{eqnarray}
\bar s^*_i  & \approx  & \frac{  a_i\left(3 \rho_i+\rho_{i+1}\right)+a_{i-1}\left(3 \rho_{i+1}-\rho_i\right) }{  8\sqrt{3}  }\, \xi  \times \\
&  & \left(  1+  \frac{
\begin{array}{c}
2 \rho_i^3 (5a_i-13a_{i-1})+\rho_{i+1} \left(2 \rho_i^2 (3a_{i-1}+5a_i)-3
  a_{i-1}a_i (3a_{i-1}+a_i)\right)+\\ 2 \rho_i \rho_{i+1}^2 (5a_{i-1}+3a_i)-
  3a_{i-1} a_i \rho_i (a_{i-1}+3a_i)+2 \rho_{i+1}^3 (5a_{i-1}-13a_i)
  \end{array}
  }{
  24a_{i-1} a_i (a_{i-1} (\rho_i+3 \rho_{i+1})+a_i (3 \rho_i+\rho_{i+1}))
  } 
  \xi ^2 + O\left[\xi\right]^4
  \right).
\end{eqnarray}
Therefore the Gaussian curvature is 
\begin{align}
\label{eqn:Ksbar}
K_{\bar s_i} & := \frac{\epsilon_{\bar s_i}}{\bar s*_i}  \approx  \frac{8 (\rho_{i-1}-\rho_{i+1})}{a_i (a_{i-1} (\rho_i+3 \rho_{i+1})+a_i (3 \rho_i+\rho_{i+1}))} \times \\
&  \left( 1 + 
\frac{
\begin{array}{c}
	a_{i-1}^2 (\rho_i+3 \rho_{i+1}) \left(54 a_i^2-5 \left(\rho_{i-1}^2+\rho_{i-1}
  	 (\rho_{i+1}-3 \rho_i)+3 \rho_i^2-3 \rho_i \rho_{i+1}+\rho_{i+1}^2\right)\right)+ \\
	 a_{i-1} a_i \left(54a_i^2 (3 \rho_i+\rho_{i+1})-5 \rho_{i-1}^2 (3 \rho_i+\rho_{i+1})+5 \rho_{i-1} \left(9
  	 \rho_i^2-\rho_{i+1}^2\right)-84 \rho_i^3+\right.\\
	 \left. 39 \rho_i^2 \rho_{i+1}+15 \rho_i \rho_{i+1}^2+10
   	\rho_{i+1}^3\right)+3 a_i^2 (\rho_i-\rho_{i+1}) \left(5 \rho_i^2+10 \rho_i \rho_{i+1}+13
 	  \rho_{i+1}^2\right)
	  \end{array}
   	}{
	36 a_{i-1} a_i^2 (a_{i-1} (\rho_i+3 \rho_{i+1})+a_i (3
  	 \rho_i+\rho_{i+1}))
	 }\, \xi ^2  
	 + O\left[\xi\right]^4\right).
   \end{align}

\subsection{The Zeroth-Order expansion term of SRF Equation for dual-edge  $\alpha_i$} 

The dual edge  SRF equation associated to $\alpha_i$ was give in Eq.~\ref{eqn:alphaequation}. We can calculate the series expansion of this equation in the limit $\xi\ll 1$.  We keep the lowest-order term by substituting the expressions we derived in the last section through Eqs.~\ref{eqn:Valphaisi}, \ref{eqn:Valphaisbari}, \ref{eqn:Ks} and \ref{eqn:Ksbar}. We find the zeroth-order term, 
\begin{equation}
\label{eqn:SRFalphazero}
\frac{2}{a_{i-1}+a_i} \left( \frac{\dot a_{i-1}+\dot a_i}{2} \right) = -
\frac{16 (a_{i-1} (\rho_i-\rho_{i+1})+a_i
   (\rho_i-\rho_{i-1})) (a_{i-1} (\rho_{i-1}+5 \rho_i+6
   \rho_{i+1})+3 a_i (3 \rho_i+\rho_{i+1}))}{3 a_{i-1}
   a_i (a_{i-1} (\rho_i+3 \rho_{i+1})+a_i (3
   \rho_i+\rho_{i+1})) (a_{i-1} (\rho_{i-1}+3 \rho_i)+a_i
   (3 \rho_i+\rho_{i+1}))}.
\end{equation}
 
We make the following substitutions, 
\begin{align}
a_{i+1} & = a_i + \zeta_1\\
a_{i-1} & = a_i - \zeta_0.
\end{align}
In the limit of  $\zeta_0 \rightarrow 0$ and $\zeta_1 \rightarrow 0$ we find Eq.~\ref{eqn:SRFalphazero} becomes,
\begin{equation}
\frac{\dot a_i}{a_i} = \frac{4 \left(\rho_{i-1}-2 \rho_i+\rho_{i+1}\right)  \left(\rho_{i-1}+14 \rho_i+9 \rho_{i+1}\right)
}{
3 a_i^2 \left(\rho_i+\rho_{i+1}) (\rho_{i-1}+6 \rho_i+\rho_{i+1}\right)}.
\end{equation}
In this limit we can also substitute
\begin{align}
\rho_{i-1}+14 \rho_i+9 \rho_{i+1} & \approx 24 \rho_i,\\
\rho_{i-1}+6 \rho_i+\rho_{i+1} & \approx 8 \rho_i,\\
\rho_i+\rho_{i+1} & \approx 2 \rho_i,
\end{align} 
yielding after some rearrangement, 
\begin{equation}
\frac{\dot a_i}{a_i} = 2\frac{\left(\frac{\rho_{i+1}-\rho_i}{a_i}\right) - \left(\frac{\rho_i-\rho_{i-1}}{a_i}\right)}{a_i}.
\end{equation}
We immediately recognize the right-hand side of this equation is the second derivative of $\rho$ with respect to $a$, and therefore  we recover the corresponding continuum RF equation for the axial edge in the limit, namely 
\begin{equation} 
\label{eqn:SRFalphaRF}
\underbrace{
\frac{\dot a_i}{a_i} = 2\frac{\left(\frac{\rho_{i+1}-\rho_i}{a_i}\right) - \left(\frac{\rho_i-\rho_{i-1}}{a_i}\right)}{a_i}
}_{\hbox{\em continuum-limit SRF equation}}  \ \  \Longrightarrow \  
\underbrace{
\frac{\dot{a}}{a}  =   2\, \frac{\rho''}{\rho}
}_{\hbox{\em continuum RF equation}}.
\end{equation}
This result proves half of Theorem 1.   We only need examine the $\sigma_i$-edge SRF equation and show it converges to the continuum RF equation in order  to complete the proof.  We accomplish this identification in the following section.

\section{The Dual $\sigma_i$-Edge SRF Equation at Point ${\cal O}_i$  and its Continuum Limit  }
\label{sec:sigmaSRF}
 
In this section we repeat the process we presented in the last section; however,  we examine the dual-edge SRF equation associated to the edge $\sigma_i\in {\cal F}^*$.  This edge in the dual cross-sectional polyhedron is dual to trapezoid $\overline{{\cal O}_i {\cal X}_i {\cal O}_{i+1} {\cal X}_{i+1}} \in {\cal F}$.  This is  illustrated in the lower right of Fig.~\ref{fig:dualareas}, and anchors the dual SRF equation in the upper left-hand side of Eq.~\ref{eq:drrf-rrf}. We also use Eqs.~\ref{eqn:s}-\ref{eqn:up} to examine the continuum limit of this equation. We use the definition of the dual-edge SRF equations as introduced recently in \cite{Miller:2013}, and consequently we find
\begin{equation}
\label{eqn:sigmaequation}
\frac{\dot \sigma_i}{\sigma_i} = -Rc_{\sigma_i }  = -\sum_{\ell_j\in\sigma^*_i} 2 \frac{\epsilon_{\ell_j}}{\ell^*_j} \left(  \frac{V_{\sigma_i \ell_j}}{V_{\sigma_i}}\right) = 
-2 \frac{\epsilon_{s_i}}{s^*_i} \left(  \frac{V_{\sigma_i s_i}}{V_{\sigma_i}}\right) 
-2 \frac{\epsilon_{s_{i+1}}}{s^*_{i+1}} \left(  \frac{V_{\sigma_i s_{i+1}}}{V_{\sigma_{i+1}}}\right) 
-2 \frac{{\epsilon}_{a_i}}{a^*_i} \left(  \frac{V_{\sigma_i a_i}}{V_{\sigma_i}}\right)-
-2 \frac{\epsilon_{\hat a_i}}{\hat a^*_i} \left(  \frac{V_{\sigma_i \hat a_i}}{V_{\sigma_i}}\right) . 
\end{equation}
Here we differentiated the axial edge $a_i = \overline{{\cal O}_i {\cal O}_{i+1}}$ from the axial edge $\hat a_i = \overline{{\cal X}_i {\cal X}_{i+1}}$. The fractional volumes, deficit angles and dual areas are not equal. As we did in the last Sec.~\ref{sec:alphaSRF}, we use the results of the Appendix and the definitions in \cite{Miller:2013} to examine each term of Eq.~\ref{eqn:sigmaequation} and to examine the Taylor series expansion to obtain the continuum limit expression. We do this in three steps by examining (1) the two circumcentric volume weighting factors on the right-hand side of the dual SRF equation, (2) the two Gaussian curvature expressions also on the right-hand side, and finally (3) the time derivative on the left-hand side of the equation. 

\subsection{The Dual Edge $\sigma_i$ to Trapezoid $\overline{{\cal O}_i {\cal X}_i {\cal O}_{i+1} {\cal X}_{i+1}}$ }

The dual edge $\sigma_i$ reaches from the circumcenter ${\cal C}_{3_{i-1}}$ of one frustum block, ${\cal F}_{i}$, to the circumcenter, ${\cal C}_{3_{i}}$, of the adjacent frustum block ${\cal F}_{i}$.   This dual edge is perpendicular to trapezoid  $\overline{ {\cal O}_i {\cal X}_i {\cal O}_{i+1} {\cal X}_{i+1}}$ common to these two frustum blocks. From \ref{subsec:dualsegments}, Eq.~\ref{eq:sigma_1}, we find
\begin{align}
\label{eqn:sigmai}
\sigma_i & = 2 \sigma_{\frac{1}{2}i} \\
& = \sqrt{\alpha_{base_i}^2 + m_{s_i\alpha_i}^2 - m_{s_i\sigma_i}^2} \\
& = \frac{   a_i^2 \bar s_i \left(s_i+s_{i+1}\right)   }{   \sqrt{4a_i^2-\left(s_i-s_{i+1}\right)^2}\sqrt{a_i^2\left(4s_i^2-\bar s_i^2\right)-s_i^2 \left(s_i-s_{i+1}\right)^2}}\\
& \approx  \frac{(\rho_i+\rho_{i+1})}{2 \sqrt{3}}\,\xi
\left( 1 +
\frac{7 (\rho_i-\rho_{i+1})^2-12a_i^2}{24 a_i^2}\,\xi^2+O\left(\xi ^4\right) \right). 
\end{align}

\subsection{Circumcentric Volume Weighting Factors}

The hybrid volume $V_{\sigma_i}$ is the volume of the hybrid polyhedron formed by the product of dual edge $\sigma_i\in {\cal F}^*$ and trapezoid $\overline{{\cal O}_i {\cal X}_i {\cal O}_{i+1} {\cal X}_{i+1}}\in {\cal F}$ and is the sum of there reduced hybrid tetrahedra,
\begin{equation}
V_{\sigma_i} = \underbrace{\frac{1}{6} \sigma_i s_i m_{\sigma_i s_i}}_{V_{\sigma_i s_i} }  +
\underbrace{\frac{1}{6} \sigma_i s_{i+1} m_{\sigma_i s_{i+1}}}_{V_{\sigma_i s_{i+1}}} + 
\underbrace{\frac{1}{6} \sigma_i a_i m_{\sigma_i a_i}}_{V_{\sigma_i a_i}} + 
\underbrace{\frac{1}{6} \sigma_i \hat a_i m_{\sigma_i \hat a_i}}_{V_{\sigma_i \hat a_i} },
\end{equation}
where the moment arms shown in the lower right-hand part of Fig.~\ref{fig:dualareas} and by Eqs.~\ref{eqn:ms1sig1}-\ref{eqn:ma1sig1}. Additionally, symmetry of the trapezoid guarantees that $m_{\sigma_i \hat a_i}= m_{\sigma_i a_i}$ so that $V_{\sigma_i a_i}=V_{\sigma_i \hat a_i}$. The fractional volumes and their series expansion in terms of $\xi$ are, 
 \begin{eqnarray}
 \label{eqn:Vsigma1s1}
 \frac{V_{\sigma_i s_i}}{V_{\sigma_i}} & =& \frac{ s_i m_{\sigma_i s_i} }{ 2\, \sigma^*_i }   = \frac{s_i\left(2 a_i^2-s_i(s_i-s_{i+1})\right)}{(s_i+s_{i+1}) \left(4 a_i^2-\left(s_i-s_{i+1}\right)^2\right)} 
 \approx  \frac{\rho_i}{2 (\rho_i+\rho_{i+1})} \left( 1+\frac{(\rho_{i+1}-\rho_i)^2}{4 a_i^2}\,\xi ^2+O\left(\xi ^4\right)\right)\\
 \label{eqn:Vsigma1s2}
 \frac{V_{\sigma_i s_{i+1}}}{V_{\sigma_i}} &  =& \frac{ s_{i+1} m_{\sigma_i s_{i+1}} }{ 2\, \sigma^*_i } = \frac{s_{i+1}\left(2 a_i^2-s_{i+1}(s_{i+1}-s_{i})\right)}{(s_i+s_{i+1}) \left(4 a_i^2-\left(s_i-s_{i+1}\right)^2\right)} 
   \approx  \frac{\rho_{i+1}}{2 (\rho_i+\rho_{i+1})} \left( 1+\frac{(\rho_i-\rho_{i+1})^2}{4
   a_i^2}\,\xi ^2+O\left(\xi ^4\right)\right),\\
    \label{eqn:Vsigma1a1}
 \frac{V_{\sigma_i a_i}}{V_{\sigma_i}} & =& \frac{ a_i m_{\sigma_i a_i} }{ 2\, \sigma^*_i } = \frac{a_i^2}{4 a_i^2-\left(s_i-s_{i+1}\right)^2}  
 \approx  \frac{1}{4} \left(1 +\frac{(\rho_i-\rho_{i+1})^2}{4 a_i^2}\,\xi ^2+O\left(\xi ^4\right)\right)\\   
 \end{eqnarray}

\subsection{Gaussian Curvature}

In \cite{Miller:2013} we defined the Gaussian curvature of edge $a_i=\overline{{\cal O}_i {\cal O}_{i+1}}$ to be the deficit angle, $\epsilon_{a_i}$, distributed uniformly over the circumcentric dual area $a^*_i$, 
\begin{equation}
K_{a_i} := \frac{\epsilon_{a_i}}{a^*_i}.
\end{equation}
The edge, $a_i$ is common to six identical isosceles frustum blocks as shown in Fig.~\ref{fig:aiedge}.
We use dihedral angle of the frustum block at edge $a_i$ given in Eq.~\ref{eqn:thetasi} to determine define the deficit angle, 
\begin{align}
\epsilon_{a_i} & = 2\pi -6 \arccos{\left( \frac{2a_i^2 \left(\bar s_i^2-2 s_i^2\right)+s_i^2
   (s_i-s_{i+1})^2}{s_i^2 \left((s_i-s_{i+1})^2-4 a_i^2\right)} \right)}\\
& \approx  
	\frac{\sqrt{3} \left(a_i^2-\left(\rho_{i+1}-\rho_i\right)^2\right)}{2a_i^2}\,\xi ^2  
		\left(1 -
		\frac{5 \left(3 a_i^4-4 a_i^2(\rho_i-\rho_{i+1})^2+(\rho_i-\rho_{i+1})^4\right)}{24  a_i^4 \left(a_i^2-\left(\rho_{i+1}-\rho_i\right)^2\right)}\, \xi ^2+O\left(\xi ^4\right)
		\right).
\end{align}
The boundary dual area $a^*_i$ is a hexagon and is shown in the upper left diagram of Fig.~\ref{fig:dualareas} and is the sum of six isosceles triangles, 
\begin{align}
a_i^* & = 6\left( \frac{1}{2}  \sigma_i m_{a_i \sigma_i} \right) 
 =  \frac{  
3 a_i^3 \bar s_i \left(s_i+s_{i+1}\right)^2
}{
2 \left(4 a_i^2-\left(s_i-s_{i+1}\right)^2\right)  \sqrt{ a_i^2\left(4s_i^2-\bar s_i^2\right)-s_i^2 \left(s_i-s_{i+1}\right)^2} 
}\\
& \approx 
\frac{ \sqrt{3}}{8} \left(\rho_i+\rho_{i+1}\right)^2\,\xi ^2  
	\left(1- \frac{2 a_i^2-5  \left(\rho_i-\rho_{i+1})^2\right)}{12 a_i^2}\,\xi ^2+
	O\left(\xi ^4\right)\right).
\end{align}
The series expansion for the Gaussian curvature of axial edge $a_i$ in powers $\xi$ is, 
\begin{align}
\label{eqn:Ka1}
K_{a_i} & := \frac{\epsilon_{a_i}}{a^*_i} \approx 
	\frac{  4 \left(a_i^2-\left(\rho_{i+1}-\rho_i\right)^2\right) }{  a_i^2(\rho_i+\rho_{i+1})^2  } 
	\left(1 +
	\frac{  5  \left(a_i^2-\left(\rho_{i+1}-\rho_i\right)^2\right) }{  24 a_i^2  }\,\xi ^2 +
	O\left(\xi ^4\right)
	\right).
\end{align}
We also need to calculate the slightly more involved expression for the Gaussian curvature, 
\begin{equation}
K_{\hat a_i} := \frac{\epsilon_{\hat a_i}}{\hat a*_i},
\end{equation}
associated with the hinge edge $\hat a_i$.  The dual area and dihedral angles now extend into the next band of triangles away from the pole ${\cal O}_i$ as shown in Fig.~\ref{fig:deltatheta}.
Using the expressions in the Appendix and the series expansion for the edges (Eq.~\ref{eqn:s}) the deficit angle for $\hat a_i$ is 
\begin{align}
\epsilon_{\hat a_i} & = 
2\pi -
2 \underbrace{
	\arccos{
		\left(
			\frac{\bar s_i \left(2 a_i^2-(s_i-s_{i+1})^2\right)}{\sqrt{4 a_i^2-(s_i-s_{i+1})^2} \sqrt{4 a_i^2
   			s_i^2-\bar s_i^2 (s_i-s_{i+1})^2}}
		\right)}
	}_{\arccos{\left(\theta_{\hat a_i\, ss}\right)}} - \\
& 2 \underbrace{
	\arccos{
		\left(
			\frac{  \bar s_i \left(2 a_i^2-(u_i-u_{i+1})^2\right)  }{  \sqrt{4a_i^2-(u_i-u_{i+1})^2} 
			\sqrt{4 a_i^2u_i^2-\bar s_i^2 (u_i-u_{i+1})^2}}
		\right)}
		}_{\arccos{\left(\theta_{\hat a_i\, uu}\right)}} - \\
&2 \underbrace{
	 \arccos{
	 	\left( 
			\frac{2 a_i^2\left(u_i^2-\bar u_i^2+{u'}_i^2\right)-{u'}_i^2(u_i-u_{i+1})^2}
			{\sqrt{4 a_i^2-(u_i-u_{i+1})^2}\sqrt{4 a_i^2 u_i^2 {u'}_i^2-{u'}_i^4(u_i-u_{i+1})^2}}
		\right)}}_{\arccos{\left(\theta_{a_i\, uu'}\right)}}\\
& \approx 
	\frac{\sqrt{3} \left(a_i^2-\left(\rho_{i+1}-\rho_i\right)^2\right)   }{  2 a_i^2 }\,  \xi ^2 
	\left(1 -
	\frac{ \left(51 a_i^4-56 a_i^2 (\rho_i-\rho_{i+1})^2+5 (\rho_i-\rho_{i+1})^4\right)  }{  24 a_i^2 \left(a_i^2-\left(\rho_{i+1}-\rho_i\right)^2\right)  }\,\xi ^2 + 
	O\left(\xi ^4\right) \right),
\end{align}
and the dual area $\hat a^*_i$ is the sum of three pairs of  kite areas whose series in $\xi$ is
\begin{equation}
\hat a^*_i   \approx  
	\frac{\sqrt{3}}{8}  \left(\rho_i+\rho_{i+1}\right)^2 \,  \xi ^2 
	\left( 1 -
	\frac{ 28 a_i^2-5 (\rho_{i+1}-\rho_i)^2}{12 a_i^2}\,\xi ^2 +
	O\left(\xi ^4\right) \right).
\end{equation}
Therefore, the Gaussian curvature at edge $\hat a_i$ is 
\begin{align}
\label{eqn:Kahat1}
K_{\hat a_i} & := \frac{\epsilon_{\hat a_i}}{\hat a^*_i}  \approx    
	-\frac{4 (\rho_i-\rho_{i+1})^2}{a_i^2 (\rho_i+\rho_{i+1})^2} 
	\left(1 -
	\frac{5 (\rho_{i+1}-\rho_i)^2  }{  24 a_i^2}\,\xi ^2 + 
	O\left(\xi ^4\right)\right).
\end{align}
Finally, the last of the four Gaussian curvatures we need for Eq.~\ref{eqn:sigmaequation} is associated with edge $s_{i+1}$. To calculate 
\begin{equation}
\label{eqn:Ks2}
K_{s_{i+1}} = \frac{\epsilon_{s_{i+1}}}{s^*_{i+1}},
\end{equation}
we need only increment each index in  Eq.~\ref{eqn:Ks} by one. 

\subsection{The Zeroth-Order Expansion Term of SRF Equation for Dual-Edge  $\sigma_i$} 

The dual edge $\sigma_i$ SRF equation was give in Eq.~\ref{eqn:sigmaequation}. We can calculate series expansion of this equation in $\xi$ and keep the lowest-order term by substituting the expressions we derived in the last section through Eqs.~\ref{eqn:Vsigma1s1}, \ref{eqn:Vsigma1s2}, \ref{eqn:Vsigma1a1} for the three circumcentric volume weighting factors, as well as the three Gaussian curvatures given in Eqs.~\ref{eqn:Ks}, \ref{eqn:Ks2}, \ref{eqn:Ka1} and \ref{eqn:Kahat1}. We find a zeroth-order term, 
\begin{eqnarray}
\label{eqn:SRFsigmazero}
\frac{\dot \rho_{i-1}+\dot \rho_i}{\rho_{i-1}+\rho_i}  & = &
	\frac{2 (\rho_i-\rho_{i+1})^2}{a_i^2 (\rho_i+\rho_{i+1})^2}
-\frac{2 (a_i+\rho_i-\rho_{i+1})
  	 (a_i-\rho_i+\rho_{i+1})}{a_i^2 (\rho_i+\rho_{i+1})^2}\\
& & -\frac{8 \left(\frac{\rho_i (a_{i-1} (\rho_i-\rho_{i+1})+a_i
   (\rho_i-\rho_{i-1}))}{a_{i-1} (a_{i-1} (\rho_{i-1}+3
   \rho_i)+a_i (3 \rho_i+\rho_{i+1}))}+\frac{\rho_{i+1} (a_i
   (\rho_{i+1}-\rho_{i+2})+a_{i+1} (\rho_{i+1}-\rho_i))}{a_{i+1} (a_i
   (\rho_i+3 \rho_{i+1})+a_{i+1} (3
   \rho_{i+1}+\rho_{i+2}))}\right)}{a_i (\rho_i+\rho_{i+1})}  
\end{eqnarray}
 We make the following substitutions, 
\begin{align}
a_{i+1} & = a_i + \zeta_1\\
a_{i-1} & = a_i - \zeta_0,
\end{align}
in the last term.  In the limit of  $\zeta_0 \rightarrow 0$ and $\zeta_1 \rightarrow 0$, we find Eq.~\ref{eqn:SRFsigmazero} becomes,
\begin{eqnarray}
\frac{\dot \rho_{i-1}+\dot \rho_i}{\rho_{i-1}+\rho_i}   & = &
	-\frac{8 \rho_i (-\rho_{i-1}+2 \rho_i-\rho_{i+1})}{a_i^2 (\rho_i+\rho_{i+1}) (\rho_{i-1}+6 \rho_i+\rho_{i+1})}-
	\frac{8 \rho_{i+1} (-\rho_i+2\rho_{i+1}-\rho_{i+2})}{a_i^2 (\rho_i+\rho_{i+1}) (\rho_i+6 \rho_{i+1}+\rho_{i+2})}-\\
& &	\frac{2\left(a_i^2-(\rho_i-\rho_{i+1})^2\right)}{a_i^2 (\rho_i+\rho_{i+1})^2}-
	\frac{2 (a_i+\rho_i-\rho_{i+1})(a_i-\rho_i+\rho_{i+1})}{a_i^2 (\rho_i+\rho_{i+1})^2}.
   \end{eqnarray}
In this limit we can also substitute
\begin{align}
\rho_{i-1}+6 \rho_i+ \rho_{i+1} & \approx 8 \rho_i,\\
\rho_{i}+6 \rho_{i+1}+\rho_{i+2} & \approx 8 \rho_{i+1},\\
\rho_i+\rho_{i+1} & \approx 2 \rho_i,
\end{align} 
in the previous equation to yield, after some rearrangement, 
\begin{equation}
\label{eqn:srfsigalmost}
\frac{\dot \rho_i}{\rho_i} = 
\frac{  
	 \frac{  
	 	\left(   \frac{   \rho_{i+1}-\rho_i}{a_i} \right)-  \left(  \frac{  \rho_i-\rho_{i-1}     }{   a_i  } \right) 
	}{ 
		a_i  
	}
}{   
	\rho_i 
} +
\frac{\left(\frac{\rho_{i+1}-\rho_i}{a_i}\right)^2}{\rho_i^2} 
-\frac{1}{\rho_i^2}.
\end{equation}
We immediately recognize the numerator of the first two terms on the right-hand side of Eq.~\ref{eqn:srfsigalmost} as the second derivative of $\rho$ with respect to $a$ and the square of the first derivative of $\rho$ with respect to $a$.  Therefore, we recover the corresponding continuum RF equation for the axial edge in the continuum limit, namely 
\begin{equation} 
\underbrace{
\frac{\dot \rho_i}{\rho_i} = \frac{  
	 \frac{  
	 	\left(   \frac{   \rho_{i+1}-\rho_i}{a_i} \right)-  \left(  \frac{  \rho_i-\rho_{i-1}     }{   a_i  } \right) 
	}{ 
		a_i  
	}
}{   
	\rho_i 
} +
\frac{\left(\frac{\rho_{i+1}-\rho_i}{a_i}\right)^2}{\rho_i^2} 
-\frac{1}{\rho_i^2}
}_{\hbox{\em continuum-limit SRF equation}} 
\ \Longrightarrow \ 
\underbrace{
\frac{\dot{\rho}}{\rho}  =  \rho''/\rho + \left(\rho'/\rho\right)^2 - 1/\rho^2
}_{\hbox{\em continuum RF equation}}.
\end{equation}
This result, together with Eq.~\ref{eqn:SRFalphaRF} and Appendix B of \cite{Alsing:2013} completes the proof of Theorem 1. 

\section{Exploring SRF as Guide for the Behavior of Hamilton's RF in Higher Dimensions}
\label{sec:explore}

We demonstrated that the continuum limit of the SRF equations yielded the Hamilton RF equations for an interesting class of warped product metrics.  Therefore, all the mathematical foundations, definitions and theorems of the continuum equations can be automatically transferred to the discrete in SRF equations.  This further reinforces the definition of the SRF equations recently forwarded as Definition 1 in \cite{Miller:2013}.  While we proved this for geometries described by the warped-product metrics in Angenent, Isenberg and  Knopf \cite{Knopf:2004,Knopf:2011}, we conjecture that the SRF equations suitably converge to the continuum Hamilton RF equations for any n-dimensional geometry for $n\ge 2$. 

We explore Ricci flow in 3 and higher dimensions because we can very well believe that the SRF equations will have an equally rich spectrum of application as does 2-dimensional  combinatorial RF applications.  We therefore are motivated to explore the discrete RF in higher dimensions so that it can be used in the analysis of topology and geometry, both numerically and analytically to bound Ricci curvature in discrete geometries and to analyze and handle higher--dimensional RF singularities \cite{LinYau:2010,Knopf:2009}.  The  topological taxonomy afforded by RF is richer when transitioning from 2 to 3--dimensions.  In particular,  the uniformization theorem says that any  2--geometry will evolve under RF to a constant curvature sphere, plane or hyperboloid, while in 3--dimensions the curvature and surface will diffuse into a connected sum of prime manifolds \cite{Thurston:1997}.  We are motivated by Alexandrov \cite{Alexandrov:1950},
\begin{quote}
{\em``The theory of polyhedra and related geometrical methods are attractive not only in their own right.  They pave the way for the general theory of surfaces.  Surely, it is not always that we may infer a theorem for curved surfaces from a theorem about polyhedra by passage to the limit.  However, the theorems about polyhedra always drive us to searching similar theorems about curved surfaces.'' A. D. Alexandrov, 1950}
\end{quote}  
We are therefore eager to explore the geometry and curvature of higher-dimensional polytopes through the SRF equations as guide for continuum analogues.

\section*{Acknowledgments} 

 WAM would like to thank the Department of Mathematics at Harvard University for their support and hospitality.  WAM and SR also acknowledge support from USAF Grant \# FA8750-11-2-0089, and support from the AFRL/RITA through the Visiting Faculty Research Program administered through the Griffiss Institute.   Any opinions, findings and conclusions or recommendations expressed in this material are those of the author(s) and do not necessarily reflect the views of the AFRL.

\begin{appendix}
\section{The Circumcentric Geometry of the Triangle-Based Frustum Block}\label{app:1}

\begin{figure}
 \includegraphics[width=14.5cm]{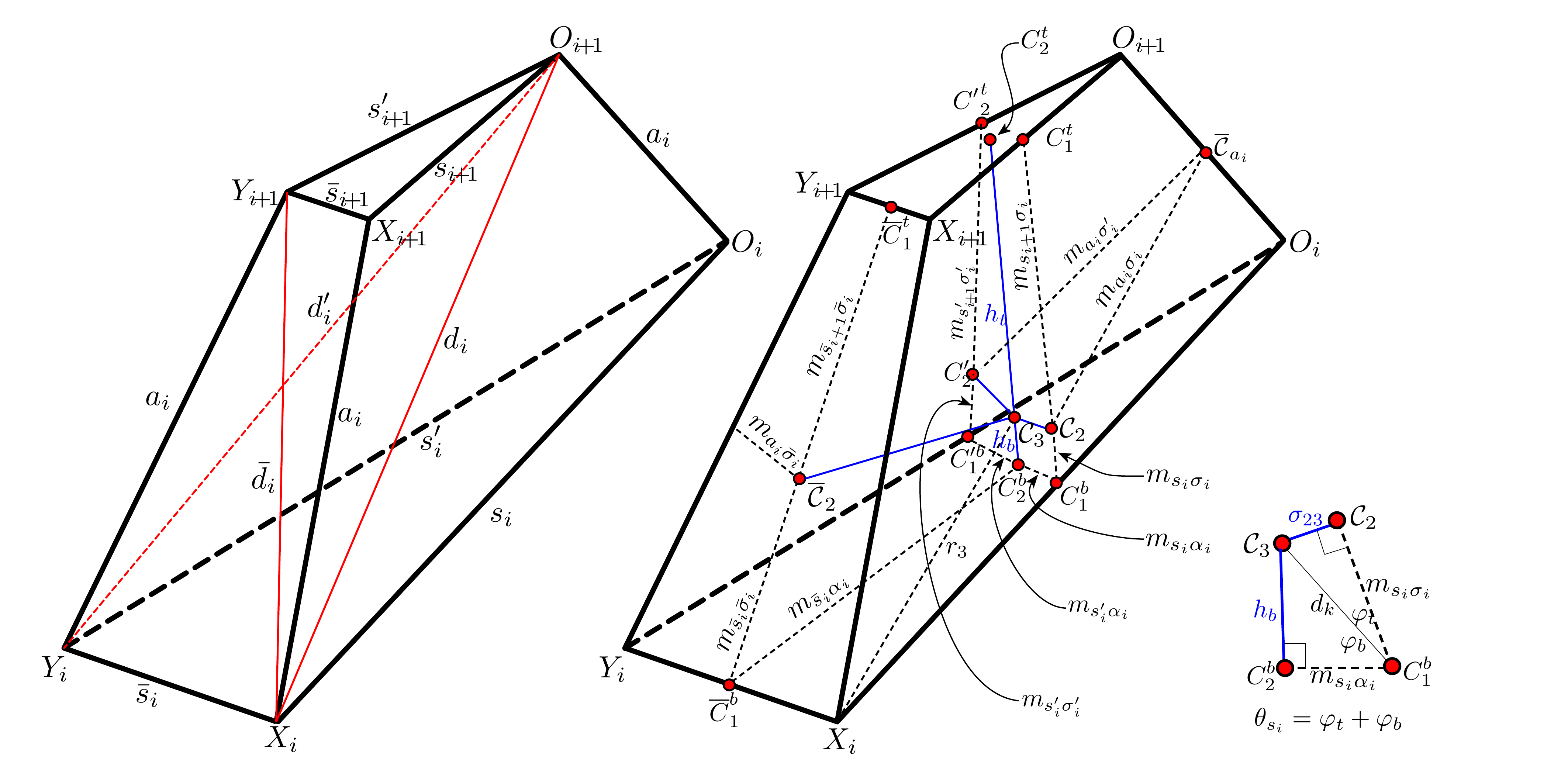}
\caption{The circumcentric-geometry of a triangle-based frustum block used to define the SRF equations}
\label{fig:frustum}
\end{figure}

The SRF equations depend on the geometry of the triangular frustum block as well as the two end-cap polyhedra. We highlight here the relevant geometric features of these frustum polyhedra.  These polyhedra are used to construct the axial-edge ($\alpha_i$) and sphere-edge ($\sigma_i$) dual SRF equations for this model.  We focus in this section on a single isosceles-based triangular frustum block as illustrated in Fig.~\ref{fig:frustum}. 
This polyhedron has nine edges.  It has three equal axial edges of length, 
\begin{equation}
\overline{O_iO_{i+1}}=\overline{X_iX_{i+1}}=\overline{Y_iY_{i+1}}=a_i,
\end{equation}
The three edges of the base triangle, $\overline{O_iX_i}=s_i$, $\overline{O_iY_i}=s'_i$,  and $\overline{X_iY_i}=\bar s_i$ can be arbitrary as long as they satisfy the triangle inequality.  We require that the top triangle is congruent, parallel and aligned with the bottom triangle.  The remaining three edges of the top triangle are $\overline{O_{i+1}X_{i+1}}=s_{i+1}$, $\overline{O_{i+1}Y_{i+1}}=s'_{i+1}$,  and $\overline{X_{i+1}Y_{i+1}}=\bar s_{i+1}$.   In this section, we assume the base triangle is larger than the top cap triangle, $s_i > s_{i+1}$, although the relevant formulas for the SRF equations will be insensitive to this choice, as we expect.  Since all three of the axial edges are equal and the triangles are congruent, the top triangle is parallel to the base triangle.  Furthermore, the circumcenter of the top triangle, ${\cal C}^t_2$, the circumcenter of the frustum, ${\cal C}_3$, and the circumcenter of the base triangle, ${\cal C}^b_2$ are collinear.  We also assume that there is no twist of the top triangle with respect to the base triangle, i.e. $s_i \parallel s_{i+1}$, i.e. the three trapezoidal faces of the frustum blocks are planar. These conditions are consistent with our model's symmetry and serve to rigidify the frustum block. 

\subsection{The nine dihedral angles of the frustum block.}

\begin{figure}
 \includegraphics[width=6 cm]{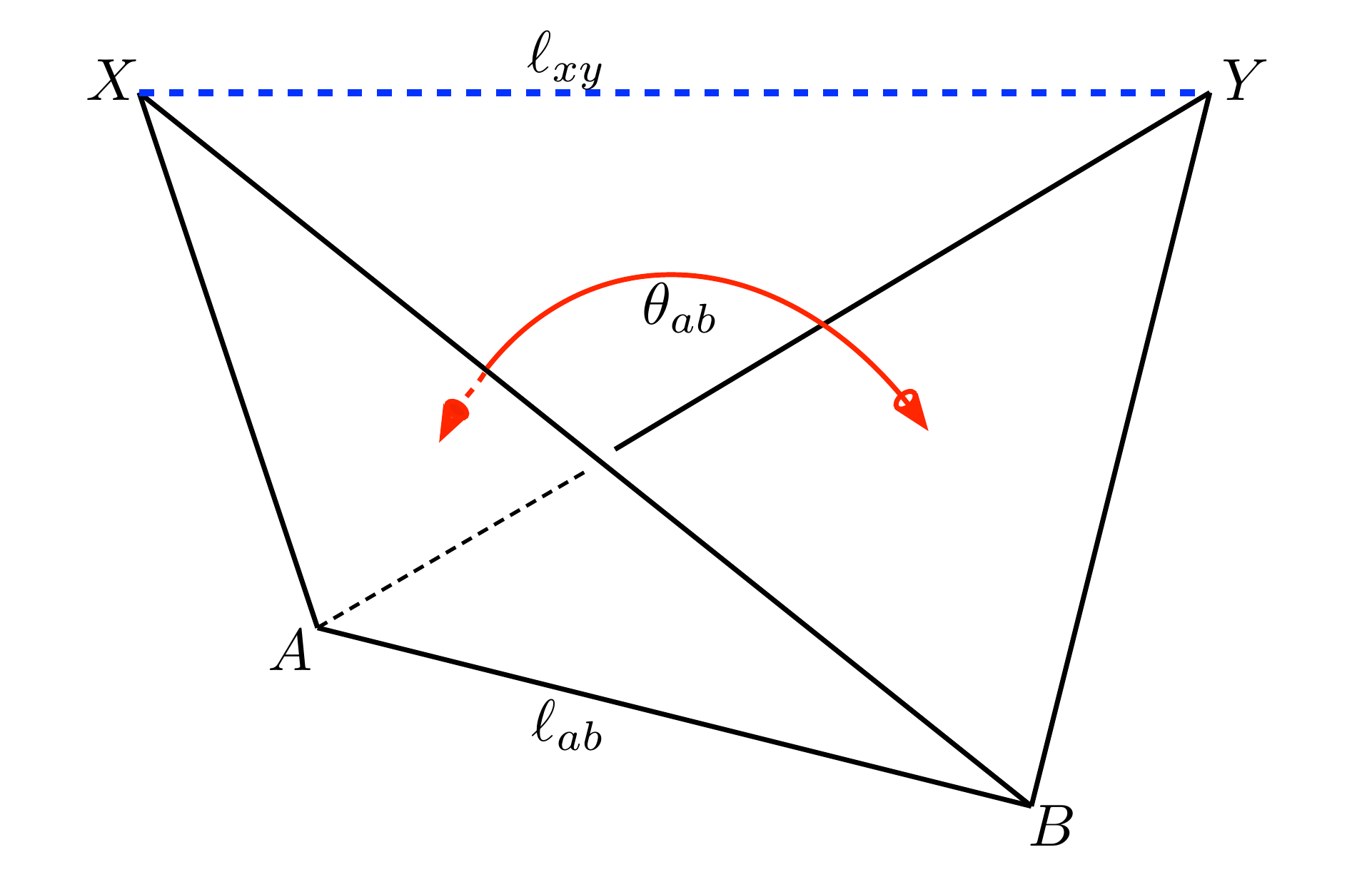}
\caption{The dihedral angle along edge $\ell_{ab}=\overline{AB}$ of tetrahedron $\{ABXY\}$.}
\label{fig:dihedral}
\end{figure}

The SRF equations are constructed, in part,  from the nine distinct dihedral angles of the frustum block.  We found it convenient to construct a diagonal for each of the three trapezoidal faces of the frustum, 
\begin{eqnarray}
d_i & = \sqrt{a_i^2 + s_i s_{i+1}}, \\
d'_i & = \sqrt{a_i^2 + s'_i s'_{i+1}}, \\
\bar d_i & = \sqrt{a_i^2 + \bar s_i \bar s_{i+1}}.
\end{eqnarray}
We can use these diagonals to subdivide the frustum block into three tetrahedra, $\{O_iO_{i+1}X_iY_i\},$ $\{Y_iY_{i+1}X_iO_{i+1}\},$  and $\{X_iX_{i+1}Y_{i+1}O_{i+1}\}.$  We can then use the usual formula for the dihedral angle of a tetrahedron to determine the nine dihedral angles of the frustum block.  The cosine of the dihedral angle for the tetrahedron shown in 
Fig.~\ref{fig:dihedral} is a function of its six edge lengths, 
\begin{equation}
\cos{\theta_{ab}} = \frac{\frac{1}{2} \left| 
\begin{array}{cc} 
 \ell_{ax}^2+\ell_{ay}^2-\ell_{xy}^2 &   \ell_{ay}^2+\ell_{bx}^2-\ell_{ab}^2-\ell_{xy}^2 \\
 \ell_{by}^2+\ell_{ax}^2-\ell_{ab}^2-\ell_{xy}^2 &  \ell_{bx}^2+\ell_{by}^2-\ell_{xy}^2
 \end{array} 
 \right|}{4 \bigtriangleup_{abx} \bigtriangleup_{aby}}. 
\end{equation} 
Here,  $\bigtriangleup_{abx}$ is the area of the triangle face $\{ABX\}$, and 
 $\bigtriangleup_{aby}$ is the area of the triangle face $\{ABY\}$.
Using this cosine formulae together with the decomposition of the frustum into three tetrahedra and the expressions for the diagonals of the three trapezoidal faces, we find the following three dihedral angles associates to the three edges of the base triangle $\bigtriangleup_i=\{O_iX_iY_i\}$: 
\begin{eqnarray}
\label{eqn:thetasi}
\theta_{s_i} & = \cos^{-1}\left(
\frac{ \left(s_{i+1}-s_i\right) \left(s^2_i-\bar s^2_i-{s'}^2_i\right)}
{4 \bigtriangleup_{i}\, \sqrt{4 a^2_i-\left(s_i-s_{i+1}\right)^2} }
 \right) \\
 \label{eqn:thetasip}
\theta_{s'_i} & =  \cos^{-1}\left(
\frac{ \left(s'_{i+1}-s'_{i}\right) \left({s'}^2_i-s^2_i-\bar s^2_i\right)}
{4 \bigtriangleup_{i}\, \sqrt{4 a^2_i-\left(\bar s_i-\bar s_{i+1}\right)^2} }
\right) \\
\label{eqn:thetasib}
\theta_{\bar s_i} & =  \cos^{-1}\left(
\frac{ \left(\bar s_{i+1}-\bar s_{i}\right) \left(\bar s^2_i-s^2_i-{s'}^2_i\right)}
{4 \bigtriangleup_{i}\, \sqrt{4 a^2_i-\left(\bar s_i-\bar s_{i+1}\right)^2} }
\right).
\end{eqnarray}
The three dihedral angles, each associated with their corresponding edge of the top cap triangle $\bigtriangleup_{i+1} =\{O_{i+1}X_{i+1}Y_{i+1}\}$ are,  
\begin{eqnarray}
\label{eqn:thetasip1}
\theta_{s_{i+1}} & = \pi - \theta_{s_i}, \\
\label{eqn:thetaip1p}
\theta_{s'_{i+1}} & = \pi - \theta_{s'_i}, \\
\label{eqn:thetaip1b}
\theta_{\bar s_{i+1}} & = \pi - \theta_{\bar s_i}.
\end{eqnarray}
The remaining three dihedral angles, $\bar \theta_{a_i}$, $\theta'_{a_i}$ and $\theta_{a_i}$ are associated with the three axial edges, $\{O_i O_{i+1}\}$, $\{X_i X_{i+1}\}$, and $\{Y_i Y_{i+1}\}$; respectively. These are similarly derived and yield, 
\begin{eqnarray}
\bar \theta_{a_i} & = &  \cos^{-1}\left(
\frac{2 a_i^2 \bar s_i^2  \left(s_i^2 + {s'}_i^2 -\bar s_i^2 \right)-s_i^2 {s'}_i^2 \left( \bar s_i-\bar s_{i+1}\right)^2}
{\bar s_i^2 s_i s'_i\, \sqrt{4 a^2_i-\left(s_i-s_{i+1}\right)^2}  \sqrt{4 a^2_i-\left(s'_i-s'_{i+1}\right)^2} }
 \right) \\
\theta'_{a_i} & = & \cos^{-1}\left(
\frac{2 a_i^2 {s'}_i^2 \left( \bar s_i^2 + s_i^2 -{s'}_i^2 \right)-s_i^2 \bar s_i^2 \left( s'_i-s'_{i+1}\right)^2}
{{s'}_i^2 s_i \bar s_i\, \sqrt{4 a^2_i-\left(s_i-s_{i+1}\right)^2}  \sqrt{4 a^2_i-\left(\bar s_i-\bar s_{i+1}\right)^2} }
\right) \\
\theta_{a_i} & = & \cos^{-1}\left(
\frac{2 a_i^2 s_i^2 \left( \bar s_i^2 + {s'}_i^2 -s_i^2 \right)-\bar s_i^2 {s'}_i^2 \left( s_i-s_{i+1}\right)^2}
{s_i^2 \bar s_i s'_i\, \sqrt{4 a^2_i-\left(s'_i-s'_{i+1}\right)^2}  \sqrt{4 a^2_i-\left(\bar s_i-\bar s_{i+1}\right)^2} }
 \right).
\end{eqnarray}

We also find it useful to  note the  three slant heights of the trapezoids, 
\begin{eqnarray}
h_{2_i }& = \overline{{\cal C}^b_1 {\cal C}^t_1} = \frac{1}{2} \sqrt{4 a_i^2 - \left( s_i-s_{i+1}\right)^2 },\\
h'_{2_i} & = \overline{{{\cal C}'}^b_1 {{\cal C}'}^t_1} = \frac{1}{2} \sqrt{4 a_i^2 - \left( s'_i-s'_{i+1}\right)^2 },\\
\bar h_{2_i} & = \overline{\bar {\cal C}^b_1 \bar {\cal C}^t_1}  =  \frac{1}{2} \sqrt{4 a_i^2 - \left( \bar s_i - \bar s_{i+1}\right)^2 },
\end{eqnarray}
height of the frustum block,
\begin{equation}\label{eq:hfi}
h_{3_i} 
= \underbrace{
\frac{
8 a_i^2 \bigtriangleup_i^2 + \bar s_i^2 {s'}_i^2 s_i \left(s_{i+1}-s_{i}\right)
}{
4 \bigtriangleup_i \sqrt{16 a_i^2 \bigtriangleup_i^2 - \bar s_i^2 {s'}_i^2 \left(s_i-s_{i+1}\right)^2}
}
}_{\alpha_{base_i}}+
\underbrace{
\frac{
8 a_i^2 \bigtriangleup_i^2 - \bar s_i^2 {s'}_i^2 s_{i+1} \left(s_{i+1}-s_i\right)
}{
4 \bigtriangleup_i \sqrt{16 a_i^2 \bigtriangleup_i^2 - \bar s_i^2 {s'}_i^2 \left(s_i-s_{i+1}\right)^2}
}
}_{\alpha_{top_i}} \\
= \frac{ 
\sqrt{16 a_i^2  \bigtriangleup^2_i- \bar s_i^2 {s'}_i^2 \left( s_1-s_{i+1}\right)^2}
}{
4 \bigtriangleup_i
},
\end{equation}
and the circumradius of the frustum block, 
\begin{equation}
r_{3_i} = 
a_i\, \sqrt{
\frac{
4 \bigtriangleup^2_i  a_i^2 +s_i s_{i+1} \bar s^2_i {s'}^2_i 
}{
16 \bigtriangleup^2_i a_i^2-\bar s^2_i {s'}^2_i \left(s_i-s_{i+1}\right)^2,
}}
\end{equation}
as labeled in Fig.~\ref{fig:frustum}.

\subsection{The twenty moment arms of the faces of the frustum block.}

Each of our SRF equations depend on the moment arms of the frustum blocks. Consider a face, $\lambda^*$, of a given frustum block. This face can be either a triangle or a trapezoid. The frustum-tiled 3-geometry, $\cal F$, that we consider here has a circumcentric dual lattice, ${\cal F}^*$.   The face $\lambda^*\in {\cal F}$ is dual to an edge, $\lambda \in \cal F^*$.  The dual edge $\lambda$ is perpendicular to face $\lambda^*$ by construction.   the moment arm ($m_{\ell \lambda}$) associated with this face and edge reaches from the middle of the edge to the circumcenter of the face.  The circumcenter of the face is the point along dual edge where $\lambda$ that intersects face $\lambda^*$.  The face will ordinarily not be the bisector of the dual edge. There are four moment arms for each of the three trapezoidal faces, and three moment arms for each of the triangular caps of the frustum block. 

The three moment arms of the base triangle as illustrated in Fig.~\ref{fig:frustum} are, 
\begin{eqnarray}
\label{eqn:msialphai}
m_{s_i \alpha_i} & = & \frac{s_i \left( \bar s^2_i + s'^2_i - s^2_i \right)}{8 \bigtriangleup_{i}},\\ 
m_{s'_i \alpha_i} & = & \frac{s'_i \left( s^2_i + \bar s^2_i - s'^2_i \right)}{8 \bigtriangleup_{i}},\\
\label{eqn:msibaralphai}
m_{\bar s_i \alpha_i} & = & \frac{\bar s_i \left( s^2_i + s'^2_i - \bar s^2_i \right)}{8 \bigtriangleup_{i}}.
\end{eqnarray} 
The three moment arms on the top triangle is obtained from these by relabeling all indices $i \longrightarrow i+1$.

The three distinct moment arms of the trapezoidal  face $\{O_iO_{i+1}X_i X_{i+1}\}$ as illustrated in Fig.~\ref{fig:frustum}, 
\begin{eqnarray}
\label{eqn:ms1sig1}
m_{s_i \sigma_i} & = & \frac{2 a^2_i + s_i \left( s_{i+1}-s_i\right)}{2 \sqrt{4 a^2_i - \left( s_i-s_{i+1}\right)^2}},\\
\label{eqn:ms1sig2}
m_{s_{i+1} \sigma_i} & = &  \frac{2 a^2_i + s_{i+1} \left( s_{i}-s_{i+1}\right)}{2 \sqrt{4 a^2_i - \left( s_i-s_{i+1}\right)^2}},\\
\label{eqn:ma1sig1}
m_{a_i \sigma_i} & = &  \frac{a_i \left( s_{i}+s_{i+1}\right)}{2 \sqrt{4 a^2_i - \left( s_i-s_{i+1}\right)^2}}.
\end{eqnarray} 
The three distinct moment arms for the trapezoidal face $\{O_iO_{i+1}Y_i Y_{i+1}\}$ can be obtained from these by the substitutions,
\begin{eqnarray}
s_i & \longrightarrow & s'_i,\\
s_{i+1} & \longrightarrow & s'_{i+1},\\
\sigma_i & \longrightarrow & \sigma'_i.
\end{eqnarray}
Similarly, the moment arms associated with face $\{X_iX_{i+1}Y_i Y_{i+1}\}$ the top triangle are obtained from Eqs.~\ref{eqn:ms1sig1}-\ref{eqn:ma1sig1} by the substitutions,
\begin{eqnarray}
s_i & \longrightarrow & \bar s_i,\\
s_{i+1} & \longrightarrow & \bar s_{i+1},\\
\sigma_i & \longrightarrow & \bar \sigma_i.
\end{eqnarray}

\subsection{The five circumcentric dual edge segments within the frustum block.}
\label{subsec:dualsegments}
This frustum block contains partial segments of five of the dual circumcentric lattice  edges,  (1) edge $\sigma_i$ dual to the trapezoid $\overline{O_i O_{i+1} X_i X_{i+1} }$, (2)  edge $\sigma'_i$ dual to face $\overline{O_i O_{i+1} Y_i Y_{i+1} }$, (3) edge $\bar \sigma_i$ dual to the trapezoid $\overline{X_i X_{i+1} Y_i Y_{i+1} }$, (4) edge $\alpha_i$ dual to the triangle $\overline{O_i X_i Y_i }$, and finally  (5) edge $\alpha_{i+1}$ dual to the triangle $\overline{O_{i+1} X_{i+1}Y_{i+1} }$.  Each of these are illustrated in Fig:~\ref{fig:dualareas}.

There are three dual $\sigma$ edges emanating from vertex ${\cal C}_3$.  The section of the circumcentric dual edges $\sigma_{\frac{1}{2}\, i}$,  $\sigma'_{\frac{1}{2}\, i}$ and $\bar \sigma_{\frac{1}{2}\, i}$ associated with the frustum block illustrated in Fig.~\ref{fig:frustum} are the line segment that starts at  the circumcenter, $C_3$ of this frustum block and terminates at the respective trapezoidal circumcenter, ${\cal C}_2$, ${\cal C}'_2$ and $\bar {\cal C}_2$; respectively.  In particular, we find that
\begin{align}
\label{eq:sigma_1}
\sigma_{\frac{1}{2} i} & = \sqrt{\alpha_{base_i}^2 + m_{s_i\alpha_i}^2 - m_{s_i\sigma_i}^2} = 
\frac{
a_i^2 \left(s_i+s_{i+1}\right) \left( \bar s_i^2 + {s'}_i^2 -s_i^2 \right)
}{
2\, \sqrt{4 a_i^2 - \left( s_i-s_{i+1}\right)^2} \sqrt{16 a_i^2 \bigtriangleup^2_i-\bar s_i^2 {s'}^2_i \left(s_i-s_{i+1}\right)^2}
} \\
\sigma'_{\frac{1}{2} i} & = \sqrt{\alpha_{base_i}^2 + m_{s'_i\alpha_i}^2 - m_{s'_i\sigma'_i}^2} =
\frac{
a_i^2 {s'}^2_i \left(s_i+s_{i+1}\right) \left( s_i^2 + \bar s_i^2 -{s'}_i^2 \right)
}{
2\, \sqrt{4 a_i^2 s_i^2 - {s'}_i^2 \left( s_i-s_{i+1}\right)^2} \sqrt{16 a_i^2 \bigtriangleup^2_i-\bar s_i^2 {s'}^2_i \left(s_i-s_{i+1}\right)^2}
} \\
\bar \sigma_{\frac{1}{2} i} & = \sqrt{\alpha_{base_i}^2 + m_{\bar s_i\alpha_i}^2 - m_{\bar s_i\bar \sigma_i}^2} =
\frac{
a_i^2  \bar s_i \left(s_i+s_{i+1}\right) \left( s_i^2 + {s'}_i^2 -\bar s_i^2 \right)
}{
2\, \sqrt{4 a_i^2 s_i^2 - \bar s_i^2 \left( s_i-s_{i+1}\right)^2} \sqrt{16 a_i^2 \bigtriangleup^2_i-\bar s_i^2 {s'}^2_i \left(s_i-s_{i+1}\right)^2}
}.
\end{align} 
Finally,  there are two segments of dual axial edges, $\alpha$, within the frustum that terminate at the circumcenter ${\cal C}_3$.  The lengths of these two segments were given in Eq.~\ref{eq:hfi}.  Their sum gives the height, $h_{3_i}$ the frustum block. One segment, $\alpha_{base_i}$, is dual to the base isosceles triangle, $\overline{{\cal O}_i {\cal X}_i {\cal Y}_i }$, 
\begin{equation} 
\label{eqn:halfalphaibase}
\alpha_{base_i} = h_b=\overline{{\cal C}^b_2 {\cal C}_3} =
\frac{
8 a_i^2 \bigtriangleup_i^2 + \bar s_i^2 {s'}_i^2 s_i \left(s_{i+1}-s_{i}\right)
}{
4 \bigtriangleup_i \sqrt{16 a_i^2 \bigtriangleup_i^2 - \bar s_i^2 {s'}_i^2 \left(s_i-s_{i+1}\right)^2}
},
\end{equation}
while,  
\begin{equation} 
\label{eqn:halfalphaitop}
\alpha_{top_i} = h_t = \overline{{\cal C}^t_2 {\cal C}_3} =
\frac{
8 a_i^2 \bigtriangleup_i^2 - \bar s_i^2 {s'}_i^2 s_{i+1} \left(s_{i+1}-s_i\right)
}{
4 \bigtriangleup_i \sqrt{16 a_i^2 \bigtriangleup_i^2 - \bar s_i^2 {s'}_i^2 \left(s_i-s_{i+1}\right)^2}
}
\end{equation}
is the other axial segment dual to triangle $\overline{{\cal O}_{i+1} {\cal X}_{i+1} {\cal Y}_{i+1} }$.

\subsection{The nine kite areas of the frustum block: the SRF equation weighting factors}

There is one kite area associated to each of the nine edges of the frustum block.  We show explicitly the kite area associated to edge $s_i$ in the far right side of Fig.~\ref{fig:frustum}.  Each kite is the sum of two triangles.  A given kite, e.g. $\kappa_{s_i}$,  is the fraction of the dual area, $s^*_i$ within the given frustum.  Since we know the moment arms and the dual-edge segments, each of the nine kites are easily expressed in terms of the nine edge lengths of the frustum block.  Using the right-most part of Fig.~\ref{fig:frustum} as guide, we find the three kite areas associated with the three edges of the base triangle,
\begin{align}
\label{eqn:kitesbase}
\kappa_{s_i}  & = \frac{1}{2} m_{s_i\sigma_i} \sigma_{\frac{1}{2} i} + \frac{1}{2} m_{s_i\alpha_i} \alpha_{base_i}, \\
\kappa_{s'_i}  & = \frac{1}{2} m_{s'_i\sigma'_i} \sigma'_{\frac{1}{2} i} + \frac{1}{2} m_{s'_i\alpha_i} \alpha_{base_i},\\
\kappa_{\bar s_i}  & = \frac{1}{2} m_{\bar s_i\sigma'_i} \bar \sigma_{\frac{1}{2} i} + \frac{1}{2} m_{\bar s_i\alpha_i} \alpha_{base_i}.
\end{align} 
The corresponding kite areas for the top triangle are, 
\begin{align}
\label{eqn:kitestop}
\kappa_{s_{i+1}}  & = \frac{1}{2} m_{s_{i+1}\sigma_i} \sigma_{\frac{1}{2} i} + \frac{1}{2} m_{s_{i+1}\alpha_i} \alpha_{top_i}, \\
\kappa_{s'_{i+1}}  & = \frac{1}{2} m_{s'_{i+1}\sigma'_i} \sigma'_{\frac{1}{2} i} + \frac{1}{2} m_{s'_{i+1}\alpha_i} \alpha_{top_i},\\
\kappa_{\bar s_{i+1}}  & = \frac{1}{2} m_{\bar s_{i+1}\sigma'_i} \bar \sigma_{\frac{1}{2} i} + \frac{1}{2} m_{\bar s_{i+1}\alpha_i} \alpha_{top_i}.
\end{align} 
Finally, there are three corresponding kite areas, $\kappa_{\bar a_i}$,  $\kappa_{a'_i}$ and $\kappa_{a_i}$ for each of the three axial edges, $\overline{{\cal O}_i{\cal O}_{i+1}}$, $\overline{{\cal X}_i{\cal X}_{i+1}}$ and $\overline{{\cal O}_i{\cal O}_{i+1}}$; respectively. We find, 
\begin{align}
\label{eqn:kitesaxial}
\kappa_{\bar a_i}  & = \frac{1}{2} m_{a_i\sigma_i} \sigma_{\frac{1}{2} i} + \frac{1}{2} m_{a_i\sigma'_i} \sigma'_{\frac{1}{2} i}, \\
\kappa_{a'_i}  & =\frac{1}{2} m_{a_i\sigma_i} \sigma_{\frac{1}{2} i} + \frac{1}{2} m_{a_i \bar \sigma_i} \bar \sigma_{\frac{1}{2} i}, \\
\kappa_{a_i}  & = \frac{1}{2} m_{a_i\sigma'_i} \sigma'_{\frac{1}{2} i} + \frac{1}{2} m_{a_i \bar \sigma_i} \bar \sigma_{\frac{1}{2} i}.
\end{align} 

Each of these kite areas can be expressed in terms of the nine edge lengths of the frustum blocks using the expressions for the moment arms and the dual edge segments above. They are rather lengthy expression, we refrain from including them here; however, they are useful in forming the dual areas used in the definition of the Gaussian curvature.. 

\end{appendix}

\end{document}